\documentclass{amsart}
\usepackage{amsfonts, amsbsy, amsmath, amssymb}
\newtheorem{thm}{Theorem}[section]
\newtheorem{lem}[thm]{Lemma}
\newtheorem{cor}[thm]{Corollary}

\newtheorem{exmp}[thm]{Example}

\numberwithin{equation}{section}

\theoremstyle{definition}

\begin{document}

\title[Explicit Evaluation of Certain Exponential Sums]{Explicit Evaluation of Certain Exponential Sums of Quadratic Functions over $\Bbb F_{p^n}$, $p$ Odd}

\author{Sandra Draper}
\address{Department of Mathematics,
University of South Florida, Tampa, FL 33620}
\email{sdraper95@yahoo.com}

\author{Xiang-dong Hou}
\address{Department of Mathematics,
University of South Florida, Tampa, FL 33620}
\email{xhou@math.usf.edu}

\keywords{Exponential sum, Gauss sum, law of quadratic reciprocity, quadratic form} 
\subjclass{11T23}
 
\begin{abstract}
Let $p$ be an odd prime and let $f(x)=\sum_{i=1}^ka_ix^{p^{\alpha_i}+1}\in\Bbb F_{p^n}[x]$, where $0\le \alpha_1<\cdots<\alpha_k$.
We consider the exponential sum $S(f,n)=\sum_{x\in\Bbb F_{p^n}}e_n(f(x))$, where $e_n(y)=e^{2\pi i\text{Tr}_n(y)/p}$, $y\in\Bbb F_{p^n}$,
$\text{Tr}_n=\text{Tr}_{\Bbb F_{p^n}/\Bbb F_p}$. There is an effective way to compute the nullity of the quadratic form $\text{Tr}_{mn}(f(x))$
for all integer $m>0$. Assuming that all such nullities are known, we find relative formulas for $S(f,mn)$ in terms of $S(f,n)$ when $\nu_p(m)
\le \min\{\nu_p(\alpha_i):1\le i\le k\}$, where $\nu_p$ is the $p$-adic order. We also find an explicit formula for $S(f,n)$ when $\nu_2(\alpha_1)=\cdots=
\nu_2(\alpha_k)<\nu_2(n)$. These results generalize those by Carlitz and by Baumert and McEliece. Parallel results with $p=2$ were obtained in a
previous paper by the second author. 
\end{abstract}

\maketitle


\section{Introduction}

Let $p$ be a prime and let $\Bbb F_{p^n}$ be the finite fields with $p^n$ elements. Denote the trace from $\Bbb F_{p^n}$ to $\Bbb F_{p^m}$ by $\text{Tr}_{n/m}$; the trace $\text{Tr}_{n/1}$ is also denoted by $\text{Tr}_n$. Let $e_n(y)=e^{2\pi i\text{Tr}_n(y)/p}$, $y\in\Bbb F_{p^n}$.
Carlitz considered the sum $\sum_{x\in\Bbb F_{p^n}}e_n(ax^{p+1}+bx)$, where $a,b\in\Bbb F_{p^n}$. He gave explicit evaluations of this sum for $p=2$
in \cite{Car79} and for odd $p$ in \cite{Car80}. The evaluation of the sum $\sum_{x\in\Bbb F_{p^n}}e_n(ax^{p^\alpha+1})$, where $a\in\Bbb F_{p^n}$ and
$\alpha\ge 0$, was also implied by the results of Baumert and McEliece \cite{Bau72}. Also see \cite {Hel76}.

More generally, we let
\begin{equation}\label{1.1}
f(x)=\sum_{i=1}^ka_ix^{p^{\alpha_i}+1}\in\Bbb F_{p^n}[x],
\end{equation}
where $0\le\alpha_1<\cdots<\alpha_k$, and consider the sum
\[
S(f,n)=\sum_{x\in\Bbb F_{p^n}}e_n(f(x)).
\]
Note that $S(f,n)$ is an exponential sum of a quadratic form on $\Bbb F_p^n$. In fact, if we identify $\Bbb F_{p^n}$ with $\Bbb F_p^n$, then 
$\text{Tr}_n(f(x))$ is a quadratic form in the coordinates of $\Bbb F_p^n$. Conversely, every quadratic form on $\Bbb F_{p^n}$ (identified with
$\Bbb F_p^n$) can be written as $\text{Tr}_n(f(x))$ for some polynomial of the form \eqref{1.1}. The sum $S(f,n)$ is completely determined by
the canonical form of the quadratic form $\text{Tr}_n(f(x))$. The converse is true except when $p=2$ and $S(f,n)=0$. In general, there is no direct
way to identify the canonical form of $\text{Tr}_n(f(x))$ using $f(x)$; an indirect way is through the computation of $S(f,n)$.   

We can also consider a seemingly more general sum 
\[
S(f+bx,n)=\sum_{x\in\Bbb F_{p^n}}e_n(f(x)+bx),
\]
where
$b\in \Bbb F_{p^n}$. However, when $p$ is odd, $S(f+bx,n)$ follows easily from $S(f,n)$, as we will see in section 3. When $p=2$, we can write $S(f+bx,n)=S(f+cx^2,n)$, where $b=c^2$ and $f+cx^2$ is of the form \eqref{1.1}.

In general, the sum $S(f,n)$ is not explicitly known. In a recent paper \cite{Hou} by the second author, Carliz's result on the sum $\sum_{x\in\Bbb F_{2^n}}e_n(ax^3+bx)$ was extended to certain polynomials of the form (1.1) in $\Bbb F_{2^n}[x]$. In the present paper, we will establish parallel
results for odd primes $p$.

When $p$ is odd, the sum $S(f,n)$ is expressed in terms of the nullity and type (definitions in the next section) of the quadratic form $\text{Tr}_n(f(x))$. The nullity of
$\text{Tr}_n(f(x))$ is rather easy to compute, as we will see in the next section. For most of the paper, 
we will treat the nullity of $\text{Tr}_n(f)$ as known and focus on the determination
of the type of $\text{Tr}_n(f(x))$. 

We fix a polynomial $f(x)\in\Bbb F_{p^n}[x]$ of the form \eqref{1.1} and consider
the sum $S(f,mn)$ as a function of $m\in\Bbb Z^+$. We try to evaluate the type of $\text{Tr}_{mn}(f(x))$ in terms of that of
$\text{Tr}_n(f(x))$; the result gives $S(f,mn)$ in terms of $S(f,n)$. For this purpose, it suffices to assume that $m$ is a prime, say $q$. The three
cases, (i) $q\ne p,2$, (ii) $q=2$, and (iii) $q=p$, require different methods. In the first two cases, we obtain explicit formulas for $S(f,qn)$ 
in terms of
$S(f,n)$. When $q=p$, we are able to express $S(f,pn)$ in terms of $S(f,n)$ only under a very restrictive condition that $\nu_p(n)<\min\{\nu_p(\alpha_i):
1\le i\le k\}$, where $\nu_p$ is the $p$-adic order function and $\alpha_1,\dots,\alpha_k$ are given in \eqref{1.1}.

The results described above are ``relative'' formulas of $S(f,mn)$ in terms of $S(f,n)$. Explicit formulas for $S(f,n)$ itself would be more
desirable. A special situation which allows explicit evaluation of $S(f,n)$ is when $\nu_2(\alpha)=\cdots=\nu_2(\alpha_k)<\nu_2(n)$. A similar
result for $p=2$ has been proved in \cite{Hou}.

The paper is organized as follows. In section 2, we express $S(f,n)$ in terms of the nullity and type of $\text{Tr}_n(f)$. We also describe a method for computing the nullity of $\text{Tr}_n(f)$. Section 3 provides a formula for $S(f+bx,n)$ in terms of $S(f,n)$. In sections 4 -- 6, we derive
relative formulas for $S(f,mn)$ in terms of $S(f,n)$ in three cases respectively: (i) $m=q^s$, where $q$ is a prime and $q\ne p,2$;
(ii) $m=2^s$; (iii) $m=p^s$, where $s\le \min\{\nu_p(\alpha_i):1\le i\le k\}-\nu_p(n)$. In section 7, we give an explicit formula for $S(f,n)$ under the assumption that $\nu_2(\alpha_1)=\cdots=\nu_2(\alpha_k)<\nu_2(n)$. As a special case of this result, we have an explicit formula for the sum
$S(ax^{p^\alpha+1})$. The two tables in the appendix contain the nullities of $\text{Tr}_n(f)$ with $p^n=3$, $\alpha_k\le 4$ and with $p^n=5$,
and $\alpha_k\le 3$. 

Throughout the paper, $p$ is always an odd prime.


\section{The Nullity and Type}

It is well known (cf. \cite[Theorem 6.9]{Jac85}) that every quadratic form $F(x_1,\dots,x_n)$ on $\Bbb F_p^n$ is congruent to a canonical form
\begin{equation}\label{2.1}
x_1^2+\cdots+x_{r-1}^2+dx_r^2,
\end{equation}
where $r\ge 0$ and $d\in\Bbb F_p^*$. (Two quartic forms $F(x_1,\dots,x_n)$ and $G(x_1,\dots,x_n)$ on $\Bbb F_p^n$ are congruent if $G(x_1,\dots,x_n)=F((x_1,\dots,x_n)A)$
for some $A\in\text{GL}(n,\Bbb F_p)$.) In \eqref{2.1}, the integer $r$ is called the rank of $F$ and $n-r$ the nullity of $F$. Let $\eta$ be the quadratic character of $\Bbb F_p$. Then
$\eta(d)\in\{\pm 1\}$ is called the type of $F$. The image of $d$ in the quotient group
$\Bbb F_p^*/(\Bbb F_p^*)^2$, where $(\Bbb F_p^*)^2=\{a^2:a\in\Bbb F_p^*\}$, is called the 
discriminant of $F$. (If $F=0$, we define the type of $F$ to $1$ and the discriminant of $F$ to be $1(\Bbb F_p^*)^2$.) The nullity of $F$ is the dimension of the $\Bbb F_p$-space
\[
\{z\in\Bbb F_p^n: F(x+z)-F(x)=\text{constant for all}\ x\in\Bbb F_p^n\}.
\]
The quadratic form $F$ is uniquely determined, up to congruence, by its rank and type.

Let $g_p$ denote the Gauss quadratic sum on $\Bbb F_p$. Recall that 
\[
g_p=
\sum_{x\in\Bbb F_p}\zeta_p^{x^2}=i^{\frac 14(p-1)^2}p^{\frac 12},
\]
where $\zeta_p=e^{2\pi i/p}$. 
 
\begin{lem}\label{L2.1}
Let $F$ be a quadratic form on $\Bbb F_p^n$ of type $t$ and rank $r$. Then
\[
\sum_{(x_1,\dots,x_n)\in\Bbb F_p^n}\zeta_p^{F(x_1,\dots,x_n)}=t\, g_p^r\, p^{n-r}.
\]
\end{lem}

\begin{proof}
We may assume $F(x_1,\dots,x_n)=x_1^2+\cdots+x_{r-1}^2+dx_r^2$, where $\eta(d)=t$. Then
\[
\sum_{(x_1,\dots,x_n)\in\Bbb F_p^n}\zeta_p^{F(x_1,\dots,x_n)}=p^{n-r}
\Bigl(\sum_{x\in\Bbb F_p}\zeta_p^{x^2}\Bigr)^{r-1}\Bigl(\sum_{x\in\Bbb F_p}\zeta_p^{dx^2}\Bigr)
=t\, g_p^r\, p^{n-r}.
\]
\end{proof}

Now let
\begin{equation}\label{2.1.0}
f(x)=\sum_{i=1}^ka_ix^{p^{\alpha_i}+1}\in\Bbb F_{p^n}[x],
\end{equation}
where $0\le \alpha_1<\cdots<\alpha_k=\alpha$ and $a_k\ne 0$. Let $l_n(f)$ and $t_n(f)$
denote the nullity and type of $\text{Tr}_n(f)$ respectively. By Lemma~\ref{L2.1},
\begin{equation}\label{2.2}
S(f,n)=t_n(f) g_p^{n-l_n(f)}p^{l_n(f)}.
\end{equation}
The nullity $l_n(f)$ is the dimension of the $\Bbb F_p$-space
\[
\{z\in\Bbb F_{p^n}:\text{Tr}_n(f(x+z)-f(x))=\text{constant for all}\ x\in\Bbb F_{p^n}\}.
\]
Note that 
\begin{equation}\label{2.2.0}
\begin{split}
&\text{Tr}_n(f(x+z)-f(x))\cr
=\,&\text{Tr}_n\Bigl[\sum_{i=1}^k a_i\bigl((x+z)^{p^{\alpha_i}+1}-x^{p^{\alpha_i}+1}\bigr)\Bigr]\cr
=\,&\text{Tr}_n\Bigl[\sum_{i=1}^k a_i\bigl(x^{p^{\alpha_i}}z+x z^{p^{\alpha_i}} +z^{p^{\alpha_i}+1}\bigr)\Bigr]\cr
=\,&\text{Tr}_n\bigl(f(z)\bigr)+\text{Tr}_n\Bigl[x\sum_{i=1}^k \bigl(a_iz^{p^{\alpha_i}}+
a_i^{p^{-\alpha_i}}z^{p^{-\alpha_i}}\bigr)\Bigr]\cr
=\,&\text{Tr}_n\bigl(f(z)\bigr)+\text{Tr}_n\Bigl[x^{p^\alpha}\sum_{i=1}^k \bigl(a_i^{p^\alpha}z^{p^{\alpha+\alpha_i}}+
a_i^{p^{\alpha-\alpha_i}}z^{p^{\alpha-\alpha_i}}\bigr)\Bigr]\cr
=\,&\text{Tr}_n\bigl(f(z)\bigr)+\text{Tr}_n\bigl[x^{p^\alpha} f^*(z)\bigr],
\end{split}
\end{equation}
where $p^{-\alpha_i}$ is a positive integer such that $p^{-\alpha_i}\cdot p^{\alpha_i}\equiv 1\pmod {p^n-1}$ and 
\begin{equation}\label{2.3}
f^*(z)=\sum_{i=1}^k\bigl(a_i^{p^\alpha}z^{p^{\alpha+\alpha_i}}+
a_i^{p^{\alpha-\alpha_i}}z^{p^{\alpha-\alpha_i}}\bigr).
\end{equation}
Therefore, $\text{Tr}_n(f(x+z)-f(x))=$ constant for all $x\in\Bbb F_{p^n}$ if and only if 
$z$ is a root of $f^*$. $f^*$ is a $p$-polynomial over $\Bbb F_{p^n}$ without repeated roots.
Thus,
\begin{equation}\label{2.4}
\begin{split}
l_n(f)\,&=\dim_{\Bbb F_p}\{z\in\Bbb F_{p^n}:f^*(z)=0\}\cr
&=\log_p\bigl|\{z\in\Bbb F_{p^n}:f^*(z)=0\}\bigr|\cr
&=\log_p\deg(f^*,\, x^{p^n}-x).
\end{split}
\end{equation}
Let 
\begin{equation}\label{2.5}
s=\min\{m:n\mid m,\ l_m(f)=2\alpha\}.
\end{equation}
Since $\deg f^*=p^{2\alpha}$, $\Bbb F_{p^s}$ is the splitting field of $f^*$ over $\Bbb F_{p^n}$. It is obvious that
\begin{equation}\label{2.6}
l_m(f)=l_{(m,s)}(f)\quad\text{for all}\ 0<m\equiv 0\pmod n.
\end{equation}
The nullity $l_m(f)$ for all $0<m\equiv 0\pmod n$ can be computed using the following 
algorithm: First use \eqref{2.4} to compute $l_{in}(f)$ for $i=1,2,\dots$ until $l_s(f)=2\alpha$.
Then $\Bbb F_{p^s}$ is the splitting field of $f^*$ and $l_m(f)$ is given by \eqref{2.6}.

\begin{exmp}\label{E2.2}
\rm
Let
\[
f(x)=x^{5^0+1}+2x^{5^1+1}+3x^{5^2+1}+4x^{5^3+1}+x^{5^4+1}\in\Bbb F_5[x].
\]
Then
\[
f^*(x)=x+4x^5+3x^{5^2}+2x^{5^3}+2x^{5^4}+2x^{5^5}+3x^{5^6}+4x^{5^7}+x^{5^8}.
\]
Using Mathematica \cite{Wol}, we find that
\begin{equation}\label{2.7}
\deg(f^*,\, x^{5^m}-x)=
\begin{cases}
5^4&\text{if}\ m=13,\cr
5^8&\text{if}\ m=26,\cr
1&\text{if}\ 1\le m\le 26,\ m\ne 13, 26.
\end{cases}
\end{equation}
Thus, the splitting field of $f^*$ over $\Bbb F_5$ is $\Bbb F_{5^{26}}$. Let $m$ be any positive integer written in the form $m=2^a13^bm^*$, where $(m^*,2\cdot 13)=1$. Then by \eqref{2.6} and \eqref{2.7},
\begin{equation}\label{2.7.0}
l_m(f)=
\begin{cases}
0&\text{if}\ b=0,\cr
4&\text{if}\ b\ge 1,\ a=0,\cr
8&\text{if}\ b\ge 1,\ a\ge 1.
\end{cases}
\end{equation}
\end{exmp}
 
In the appendix, we give the values of $l_m(f)$, where $f$ is of the form \eqref{2.1.0} with
$p^n=3$, $\alpha\le 4$ and $p^n=5$, $\alpha\le 3$.


\section{From $S(f,n)$ to $S(f+bx,n)$}

Let $f(x)\in\Bbb F_{p^n}[x]$ be as in \eqref{2.1.0}.

\begin{lem}\label{L3.1}
Let $b\in\Bbb F_{p^n}$. Then
\[
\begin{split}
&S(f,n)\overline{S(f+bx,n)}\cr
=\,&
\begin{cases}
p^{n+l_n(f)}e_n(f(x_0))&\text{if $f^*(x)=b^{p^\alpha}$ has a solution $x_0\in\Bbb F_{p^n}$}\cr
0&\text{otherwise}.
\end{cases}
\end{split}
\]
\end{lem}

\begin{proof}
We have
\[
\begin{split}
&S(f,n)\overline{S(f+bx,n)}\cr
=\,&\sum_{x,y\in\Bbb F_{p^n}}e_n\bigl(f(x)-f(y)-by\bigr)\cr
=\,&\sum_{x,y\in\Bbb F_{p^n}}e_n\bigl(f(x+y)-f(y)-by\bigr)\cr
=\,&\sum_{x,y\in\Bbb F_{p^n}}e_n\bigl(f(x)+yf^*(x)^{p^{-\alpha}}-by\bigr)\kern1cm
\text{(by \eqref{2.2.0})}\cr
=\,&\sum_{x\in\Bbb F_{p^n}}e_n\bigl(f(x)\bigr)
\sum_{y\in\Bbb F_{p^n}}e_n\Bigl(y\bigl(f^*(x)^{p^{-\alpha}}-y\bigr)\Bigr)\cr
=\,& p^n\sum_{\substack{x\in\Bbb F_{p^n}\cr f^*(x)=b^{p^\alpha}}} e_n\bigl(f(x)\bigr).
\end{split}
\]
If $f^*(x)=b^{p^\alpha}$ has no solution in $\Bbb F_{p^n}$, $S(f,n)\overline{S(f+bx,n)}=0$.
If $f^*(x)=b^{p^\alpha}$ has a solution $x_0\in \Bbb F_{p^n}$, then the solution
set of $f^*(x)=b^{p^\alpha}$ in $\Bbb F_{p^n}$ is the $l_n(f)$-dimensional affine
subspace $x_0+\{z\in\Bbb F_{p^n}:f^*(z)=0\}$ and $\text{Tr}_n(f(x))$ is a constant on this
affine subspace. Therefore,
\[
S(f,n)\overline{S(f+bx,n)}=p^{n+l_n(f)}e_n(f(x_0)).
\]
\end{proof}

\begin{cor}\label{C3.2}
We have
\[
S(f+bx,n)=
\begin{cases}
\overline{e_n(f(x_0))}S(f,n)&\text{if $f^*(x)=b^{p^\alpha}$ has a solution $x_0\in\Bbb F_{p^n}$},\cr
0&\text{otherwise}.
\end{cases}
\]
\end{cor}

\begin{proof}
Since $S(f,n)\ne 0$, by Lemma~\ref{L3.1}, we only have to consider the case when
$f^*(x)=b^{p^\alpha}$ has a solution $x_0\in\Bbb F_{p^n}$. We have
\[
\begin{split}
S(f+bx,n)\,&=\frac{p^{n+l_n(f)}}{\overline{S(f,n)}}\,\overline{e_n(f(x_0))}\cr
&=\frac{p^{n+l_n(f)}}{|S(f,n)|^2}\,\overline{e_n(f(x_0))}\,S(f,n)\cr
&= \overline{e_n(f(x_0))}\,S(f,n)\kern 1cm\text{(by \eqref{2.2})}.
\end{split}
\]
\end{proof}


\section{From $S(f,n)$ to $S(f, q^sn)$, $q\ne 2,\, p$}

In this section, we assume that $q$ is a odd prime with $q\ne p$ and $s\ge 0$.

\begin{lem}\label{L4.1}
In the ring $\Bbb Z[\zeta_p]$, we have
\begin{equation}\label{4.1}
S(f,q^sn)\equiv \sum_{x\in\Bbb F_{p^n}}e_{q^sn}\bigl(f(x)\bigr)\pmod q;
\end{equation}
equivalently,
\begin{equation}\label{4.2}
t_{q^sn}(f)g_p^{q^sn-l_{q^sn}(f)}p^{l_{q^sn}(f)}\equiv\Bigl(\frac qp\Bigr)^{s(n-l_n(f))}
t_n(f) g_p^{n-l_n(f)}p^{l_n(f)}\pmod q,
\end{equation}
where $\bigl(\frac qp\bigr)$ is the Legendre symbol (\cite{Ire82}).
\end{lem}

\begin{proof}
Write
\[
S(f,q^sn)=\sum_{x\in\Bbb F_{p^{q^sn}}\setminus\Bbb F_{p^n}}e_{q^sn}\bigl(f(x)\bigr)+
\sum_{x\in\Bbb F_{p^n}}e_{q^sn}\bigl(f(x)\bigr).
\]
In the above,
\[
\sum_{x\in\Bbb F_{p^{q^sn}}\setminus\Bbb F_{p^n}}e_{q^sn}\bigl(f(x)\bigr)=\sum_{i=0}^{p-1}
|T(i)|\zeta_p^i,
\]
where $T(i)=\{x\in\Bbb F_{p^{q^sn}}\setminus\Bbb F_{p^n}:\text{Tr}_{q^sn}(f(x))=i\}$.
$T(i)$ is a union of $\text{Aut}(\Bbb F_{p^{q^sn}}/\Bbb F_{p^n})$-orbits of cardinality
$>1$. Since the cardinality of every $\text{Aut}(\Bbb F_{p^{q^sn}}/\Bbb F_{p^n})$-orbit is
a power of $q$, we have $|T(i)|\equiv 0\pmod q$ for all $0\le i\le p-1$.
Hence,
\[
\sum_{x\in\Bbb F_{p^{q^sn}}\setminus\Bbb F_{p^n}}e_{q^sn}\bigl(f(x)\bigr)\equiv 0\pmod q.
\]
So, \eqref{4.1} is proved.

To see the equivalence between \eqref{4.1} and \eqref{4.2}, we compute both sides of \eqref{4.1}. By \eqref{2.2},
\[
S(f,q^sn)=t_{q^sn}(f)\,g_p^{q^sn-l_{q^sn}(f)}p^{l_{q^sn}(f)}.
\]
Note that the quadratic forms $q^s\text{Tr}_n(f(x))$ and $\text{Tr}_n(f(x))$ on $\Bbb F_{p^n}$
have the same rank but the type of $q^s\text{Tr}_n(f(x))$ equals $\bigl(\frac qp\bigr)^{s(n-l_n(f))}$ times the type of $\text{Tr}_n(f(x))$. Thus,
\[
\sum_{x\in\Bbb F_{p^n}}e_{q^sn}\bigl(f(x)\bigr)=\sum_{x\in\Bbb F_{p^n}}\zeta_p^{q^s\text{Tr}_n(f(x))}=\bigl(\frac qp\bigr)^{s(n-l_n(f))} t_n(f)
g_p^{n-l_n(f)}p^{l_n(f)}.
\]
Now, the equivalence between \eqref{4.1} and \eqref{4.2} is clear.
\end{proof}

\begin{lem}\label{L4.2}
We have
\begin{equation}\label{4.3}
p^{l_{q^sn}(f)-l_n(f)}\equiv 1\pmod q
\end{equation}
and
\begin{equation}\label{4.4}
\frac12\bigl[l_{q^sn}(f)-l_n(f)\bigr]\in\Bbb Z.
\end{equation}
\end{lem}

\begin{proof}
Let $X=\{x\in\Bbb F_{p^{q^sn}}\setminus\Bbb F_{p^n}:f^*(x)=0\}$. Then $X$ is a union of 
$\text{Aut}(\Bbb F_{p^{q^sn}}/\Bbb F_{p^n})$-orbits of cardinality $>1$. Thus,
\begin{equation}\label{4.5}
p^{l_{q^sn}(f)}-p^{l_n(f)}=|X|\equiv 0\pmod q
\end{equation}
and \eqref{4.3} is proved.

Using \eqref{4.3}, we can simplify \eqref{4.2} as
\begin{equation}\label{4.6}
t_{q^sn}(f)g_p^{q^sn-l_{q^sn}(f)}\equiv\Bigl(\frac qp\Bigr)^{s(n-l_n(f))}t_n(f)
g_p^{n-l_n(f)}\pmod q.
\end{equation}
Assume to the contrary of \eqref{4.4} that $l_{q^sn}(f)-l_n(f)$ is odd. Then exactly one of 
$q^sn-l_{q^sn}(f)$ and $n-l_n(f)$ is odd. Also note that $g_p^2=\pm p$. Thus \eqref{4.6}
gives
\begin{equation}\label{4.7}
g_pp^a\equiv b\pmod q
\end{equation}
for some $a,b\in\Bbb Z$, $a\ge 0$. Let $\sigma\in\text{Aut}(\Bbb Q(g_p)/\Bbb Q)$ such that
$\sigma(g_p)=-g_p$. Since $\Bbb Q(\zeta_p)/\Bbb Q$ is Galois, we can extend $\sigma$ to $\tau\in
\text{Aut}(\Bbb Q(\zeta_p)/\Bbb Q)$. Apply $\tau$ to \eqref{4.7}. We have
\[
-g_pp^a\equiv b\pmod q.
\]
It follows that $2g_pp^a\equiv 0\pmod q$, which is a contradiction.
\end{proof}

Let $o_q(p)$ denote the multiplicative order of $p$ in $\Bbb Z/q\Bbb Z$. Then by \eqref{4.3},
$\frac 1{o_q(p)}\bigl(l_{q^sn}(f)-l_n(f)\bigr)\in\Bbb Z$.

\begin{thm}\label{T4.3}
We have
\begin{equation}\label{4.7.0}
t_{q^sn}(f)=t_n(f)\Bigl(\frac qp\Bigr)^{sl_n(f)}(-1)^{\frac 14(p-1)[l_{q^sn}(f)-l_n(f)]+\frac 1{o_q(p)}[l_{q^sn}(f)-l_n(f)]}.
\end{equation}
\end{thm}

\begin{proof}
By \eqref{4.6},
\begin{equation}\label{4.8}
t_{q^sn}(f)g_p^{n(q^s-1)}\equiv\Bigl(\frac qp\Bigr)^{s(n-l_n(f))}t_n(f)g_p^{l_{q^sn}(f)-l_n(f)}\pmod q.
\end{equation}
Since $g_p^2=(-1)^{\frac {p-1}2}p$, we have
\begin{equation}\label{4.9}
\begin{split}
g_p^{n(q^s-1)}&\,=\Bigl[(-1)^{\frac {p-1}2}p\Bigr]^{\frac 12 n(q^s-1)}=(-1)^{\frac 14n(p-1)(q^s-1)}p^{\frac 12(q-1)n(1+q+\cdots+q^{s-1})}\cr
&\equiv (-1)^{\frac 14sn(p-1)(q-1)}\Bigl(\frac pq\Bigr)^{sn} \pmod q.
\end{split}
\end{equation}
(Note that $p^{\frac 12(q-1)}\equiv\bigl(\frac pq\bigr)\pmod q$ and that $q^s-1\equiv s(q-1)\pmod 4$.) Also,
\[
\begin{split}
g_p^{l_{q^sn}(f)-l_n(f)}\,&=\Bigl[(-1)^{\frac{p-1}2}p\Bigr]^{\frac 12(l_{q^sn}(f)-l_n(f))}\cr
&=(-1)^{\frac 14(p-1)(l_{q^sn}(f)-l_n(f))}p^{\frac 12(l_{q^sn}(f)-l_n(f))}.
\end{split}
\]
Since $\frac 12(l_{q^sn}(f)-l_n(f))\in\Bbb Z$ and $o_q(p)\mid (l_{q^sn}(f)-l_n(f))$ (Lemma~\ref{L4.2}), we have
\[
p^{\frac 12(l_{q^sn}(f)-l_n(f))}\equiv (-1)^{\frac 1{o_q(p)}(l_{q^sn}(f)-l_n(f))}\pmod q.
\]
Thus,
\begin{equation}\label{4.10}
g_p^{l_{q^sn}(f)-l_n(f)}\equiv (-1)^{\frac 14(p-1)(l_{q^sn}(f)-l_n(f))+\frac 1{o_q(p)}(l_{q^sn}(f)-l_n(f))} \pmod q.
\end{equation}
Using \eqref{4.9} and \eqref{4.10} in \eqref{4.8}, we have
\[
\begin{split}
t_{q^sn}(f)\equiv\,& t_n(f)\Bigl(\frac qp\Bigr)^{s(n-l_n(f))}\Bigl(\frac pq\Bigr)^{sn}\cr 
&\cdot (-1)^{\frac 14(p-1)(sn(q-1)+l_{q^sn}(f)-l_n(f))+\frac 1{o_q(p)}(l_{q^sn}(f)-l_n(f))} \pmod q.
\end{split}
\]
In the above, both sides are $\pm 1$; hence the two sides are equal. Using the law of quadratic reciprocity
\[
\Bigl(\frac pq\Bigr)\Bigl(\frac qp\Bigr)=(-1)^{\frac 14(p-1)(q-1)},
\]
we have \eqref{4.7.0}.
\end{proof}

\begin{exmp}\label{E4.4}
\rm
(Example~\ref{E2.2} continued)
Let 
\[
f(x)=x^{5^0+1}+2x^{5^1+1}+3x^{5^2+1}+4x^{5^3+1}+x^{5^4+1}\in\Bbb F_5[x]
\]
be the polynomial considered in Example~\ref{E2.2}. The nullity $l_m(f)$ is given by \eqref{2.7.0} for all $m>0$. Since
\[
S(f,1)=\sum_{x\in\Bbb F_5}e_1(f(x))=\sum_{x\in\Bbb F_5}e_1(x^2)=g_5,
\]
we have
\[
t_1(f)=1.
\]
Now let $m$ be odd and not divisible by $5$. Write $m=13^bm^*$, where $(m^*,2\cdot 5\cdot 13)=1$. Note that $o_{13}(5)=4$. Then it follows from \eqref{4.7.0} that 
\[
t_m(f)=
\begin{cases}
1&\text{if}\ b=0,\cr
-1&\text{if}\ b\ge 1.
\end{cases}
\]
In the next section, we will revisit this example and determine $t_m(f)$ for all even $m$ not divisible by $5$.
\end{exmp}


\section{From $S(f,n)$ to $S(f,2^sn)$, $s>0$}

We have
\[
S(f,2^sn)=\sum_{x\in\Bbb F_{p^{2^sn}}\setminus\Bbb F_{p^{2n}}}e_{2^sn}(f(x))+\sum_{x\in\Bbb F_{p^{2n}}}e_{2^sn}(f(x)).
\]
Here,
\[
\sum_{x\in\Bbb F_{p^{2^sn}}\setminus\Bbb F_{p^{2n}}}e_{2^sn}(f(x))=\sum_{i=0}^{p-1}|T(i)|\zeta_p^i,
\]
where
\[
T(i)=\{x\in\Bbb F_{p^{2^sn}}\setminus\Bbb F_{p^{2n}}: \text{Tr}_{2^sn}(f(x))=i\}.
\]
We claim that $T(i)$ is a union of $\text{Aut}(\Bbb F_{p^{2^sn}}/\Bbb F_{p^n})$-orbits of cardinality divisible by $4$. Note that 
$\text{Aut}(\Bbb F_{p^{2^sn}}/\Bbb F_{p^n})$ is cyclic of order $2^s$. If $x\in\Bbb F_{p^{2^sn}}\setminus\Bbb F_{p^{2n}}$,
its stabilizer in $\text{Aut}(\Bbb F_{p^{2^sn}}/\Bbb F_{p^n})$ does not contain $\text{Aut}(\Bbb F_{p^{2^sn}}/\Bbb F_{p^{2n}})$ hence 
must be properly contained in $\text{Aut}(\Bbb F_{p^{2^sn}}/\Bbb F_{p^{2n}})$. So, the $\text{Aut}(\Bbb F_{p^{2^sn}}/\Bbb F_{p^n})$-orbit
of $x$ has a cardinality divisible by $4$; hence the claim is proved.

Therefore, we have
\begin{equation}\label{5.1}
S(f,2^sn)\equiv\sum_{x\in\Bbb F_{p^{2n}}}e_{2^sn}(f(x))\pmod 4.
\end{equation}
Partition $\Bbb F_{p^{2n}}$ as
\[
\Bbb F_{p^{2n}}=\Bbb F_{p^n}\overset{\centerdot}\cup A\overset{\centerdot}\cup B,
\]
where
\[
\begin{split}
A\,&=\{x\in\Bbb F_{p^{2n}}^*:x^{p^n-1}=-1\},\cr
B\,&=\{x\in\Bbb F_{p^{2n}}^*:x^{p^n-1}\ne\pm 1\}.
\end{split}
\]
Note that $B$ can be further partitioned into four-element subsets of the form $\{\pm x, \pm x^{p^n}\}$. Moreover, $e_{2^sn}(f(x))$ is constant on $\{\pm x, \pm x^{p^n}\}$. Thus, 
\begin{equation}\label{5.2}
\sum_{x\in B}e_{2^sn}(f(x))\equiv 0\pmod 4.
\end{equation}
Let $\beta$ be a nonsquare of $\Bbb F_{p^n}$ and let $x_0\in\Bbb F_{p^{2n}}$ such that $x_0^2=\beta$. Then $x_0^{p^n-1}=-1$ and $A=x_0\Bbb F_{p^n}^*$. So, 
\begin{equation}\label{5.3}
\begin{split}
\sum_{x\in A}e_{2^sn}(f(x))\,&=\sum_{x\in \Bbb F_{p^n}^*}e_{2^sn}\Bigl(\sum_{i=1}^ka_ix_0^{p^{\alpha_i}+1}x^{p^{\alpha_i}+1}\Bigr)\cr
&=\sum_{x\in \Bbb F_{p^n}^*}e_{2^sn}\Bigl(\sum_{i=1}^ka_i\beta^{\frac 12(p^{\alpha_i}+1)}x^{p^{\alpha_i}+1}\Bigr)\cr
&=\sum_{x\in \Bbb F_{p^n}}e_{2^sn}\bigl(\tilde f(x)\bigr)-1,
\end{split}
\end{equation}
where
\begin{equation}\label{5.3.0}
\tilde f(x)=\sum_{i=1}^ka_i\beta^{\frac 12(p^{\alpha_i}+1)}x^{p^{\alpha_i}+1}\in \Bbb F_{p^n}[x].
\end{equation}
Note that when $\beta$ is fixed, $\tilde{\tilde f}(x)=f(\beta x)$.
By \eqref{5.1} -- \eqref{5.3}, we have
\[
\begin{split}
S(f,2^sn)\,&\equiv\sum_{x\in\Bbb F_{p^n}}e_{2^sn}(f(x))+\sum_{x\in A}e_{2^sn}(f(x))\pmod 4\cr
&=\sum_{x\in\Bbb F_{p^n}}e_{2^sn}(f(x))+\sum_{x\in\Bbb F_{p^n}}e_{2^sn}(\tilde f(x))-1\cr
&=\Bigl(\frac 2p\Bigr)^{s(n-l_n(f))}S(f,n)+\Bigl(\frac 2p\Bigr)^{s(n-l_n(\tilde f))}S(\tilde f,n)-1,
\end{split}
\]
i.e.,
\[
\begin{split}
&t_{2^sn}(f)g_p^{2^sn-l_{2^sn}(f)}p^{l_{2^sn}(f)}\cr
\equiv\,&\Bigl(\frac 2p\Bigr)^{s(n-l_n(f))} t_n(f) g_p^{n-l_n(f)}p^{l_n(f)}+\Bigl(\frac 2p\Bigr)^{s(n-l_n(\tilde f))} t_n(\tilde f) g_p^{n-l_n(\tilde f)}p^{l_n(\tilde f)}\pmod 4.
\end{split}
\]
Since $\bigl(\frac 2p\bigr)=(-1)^{\frac 18(p^2-1)}$ and $p\equiv (-1)^{\frac 12(p-1)}\pmod 4$, the above can be written as 
\begin{equation}\label{5.4}
\begin{split}
&t_{2^sn}(f)(-1)^{\frac 12(p-1)l_{2^sn}(f)}g_p^{2^sn-l_{2^sn}(f)}\cr
\equiv\,&t_n(f)(-1)^{\frac 18(p^2-1)s(n-l_n(f))+\frac 12(p-1)l_n(f)}g_p^{n-l_n(f)}\cr
&+t_n(\tilde f)(-1)^{\frac 18(p^2-1)s(n-l_n(\tilde f))+\frac 12(p-1)l_n(\tilde f)}g_p^{n-l_n(\tilde f)}-1 \pmod 4.
\end{split}
\end{equation}

\begin{thm}\label{T5.1}
Let $f$ be given by \eqref{2.1.0} and $\tilde f$ by \eqref{5.3.0} and let $s>0$. Then $l_n(f)+l_n(\tilde f)+l_{2^sn}(f)$ is even.
Moreover,
\begin{equation}\label{5.5}
t_{2^sn}(f)=
\begin{cases}
t_n(f)t_n(\tilde f)&\text{if}\ l_n(f)\equiv l_n(\tilde f)\pmod 2,\cr
t_n(f)t_n(\tilde f)(-1)^{\frac 18(p^2-1)s}&\text{if}\ l_n(f)\not\equiv l_n(\tilde f)\pmod 2.
\end{cases}
\end{equation}
\end{thm}

\begin{proof}
We first show that $l_n(f)+l_n(\tilde f)+l_{2^sn}(f)$ is even. If, to the contrary, $l_n(f)+l_n(\tilde f)+l_{2^sn}(f)$ is odd, then exactly one or
three of $2^sn-l_{2^sn}(f)$, $n-l_n(f)$ and $n-l_n(\tilde f)$ are odd. Note that $g_p^2=(-1)^{\frac 14(p-1)^2}p\equiv 1\pmod 4$. Thus, \eqref{5.4} gives
\[
g_p\equiv u\pmod 4\qquad\text{for some}\ u\in\Bbb Z
\]
or
\[
-1\equiv vg_p\pmod 4\qquad\text{for some}\ v\in\Bbb Z,
\]
depending on whether one or three of $2^sn-l_{2^sn}(f)$, $n-l_n(f)$ and $n-l_n(\tilde f)$ are odd. Let $\tau\in\text{Aut}(\Bbb Q(\zeta_p)/\Bbb Q)$ such
that $\tau(g_p)=-g_p$ and apply $\tau$ to the above. In the first case, we have
\[
\begin{cases}
g_p\equiv u\pmod 4,\cr
-g_p\equiv u\pmod 4.
\end{cases}
\]
So, $2g_p\equiv 0\pmod 4$, which is a contradiction. In the second case, we have 
\[
\begin{cases}
-1\equiv vg_p\pmod 4,\cr
-1\equiv -vg_p\pmod 4.
\end{cases}
\]
So, $-2\equiv 0\pmod 4$, which is also a contradiction. Thus, we have proved that $l_n(f)+l_n(\tilde f)+l_{2^sn}(f)$ is even.

To prove \eqref{5.5}, we first assume $l_n(f)\equiv l_n(\tilde f)\pmod 2$. Then $l_{2^sn}(f)$ is even. By \eqref{5.4}, we have
\begin{equation}\label{5.6}
t_{2^sn}(f)\equiv \bigl(t_n(f)+t_n(\tilde f)\bigr)\delta-1\pmod 4,
\end{equation}
where $\delta\in\{\pm 1, \pm g_p\}$. Since $\frac 12(g_p-1)$ is integral over $\Bbb Q$, $g_p\equiv 1\pmod 2$. Thus, $2\delta\equiv 2\pmod 4$. 
If $t_n(f)+t_n(\tilde f)=0$, \eqref{5.6} gives $t_{2^sn}(f)=-1$. If $t_n(f)+t_n(\tilde f)=\pm 2$, \eqref{5.6} gives
\[
t_{2^sn}(f)\equiv 2\delta-1\equiv 1\pmod 4,
\]
i.e.,
$t_{2^sn}(f)=1$. To sum up, we have 
\[
t_{2^sn}(f)=\left\{
\begin{array}{ll}
-1&\text{if}\ t_n(f)+t_n(\tilde f)=0\cr
1&\text{if}\ t_n(f)+t_n(\tilde f)=\pm 2
\end{array}\right\}
=t_n(f)t_n(\tilde f).
\]

Now assume $l_n(f)\not\equiv l_n(\tilde f)\pmod 2$. Then $l_{2^sn}(f)$ is odd. Without loss of generality, assume $n-l_n(f)$ is odd and $n-l_n(\tilde f)$ is even. Then by \eqref{5.4},
\[
\begin{split}
&t_{2^sn}(f)(-1)^{\frac 12(p-1)}g_p\cr
\equiv\,& t_n(f)(-1)^{\frac 18(p^2-1)s+\frac 12(p-1)l_n(f)}g_p+t_n(\tilde f)(-1)^{\frac 12(p-1)l_n(\tilde f)}-1\pmod 4\cr
=\,&t_n(f)(-1)^{\frac 18(p^2-1)s+\frac 12(p-1)(l_n(\tilde f)+1)}g_p+t_n(\tilde f)(-1)^{\frac 12(p-1)l_n(\tilde f)}-1,
\end{split}
\]
i.e.,
\[
\begin{split}
&(-1)^{\frac 12(p-1)}\Bigl[t_{2^sn}(f)-t_n(f)(-1)^{\frac 18(p^2-1)s+\frac 12(p-1)l_n(\tilde f)}\Bigr]g_p\cr
\equiv\,& t_n(\tilde f)
(-1)^{\frac 12(p-1)l_n(\tilde f)}-1\pmod 4.
\end{split}
\]
So, 
\[
\begin{split}
t_{2^sn}(f)\,&=
\begin{cases}
\displaystyle t_n(f)(-1)^{\frac 18(p^2-1)s+\frac 12(p-1)l_n(\tilde f)}&\text{if}\ t_n(\tilde f)(-1)^{\frac 12(p-1)l_n(\tilde f)}=1,
\vrule height0pt width0pt depth 8pt\cr
\displaystyle -t_n(f)(-1)^{\frac 18(p^2-1)s+\frac 12(p-1)l_n(\tilde f)}&\text{if}\ t_n(\tilde f)(-1)^{\frac 12(p-1)l_n(\tilde f)}=-1
\end{cases}\cr
&=t_n(f)t_n(\tilde f)(-1)^{\frac 18(p^2-1)s}.
\end{split}
\]
This completes the proof of the theorem.
\end{proof}

\begin{exmp}\label{E5.2}
\rm
(Example~\ref{E4.4} revisited)
Recall that 
\[
f(x)=x^{5^0+1}+2x^{5^1+1}+3x^{5^2+1}+4x^{5^3+1}+x^{5^4+1}\in\Bbb F_5[x],
\]
and $t_1(f)=1$, $l_1(f)=0$. Choose a nonsquare $\beta=2\in\Bbb F_5$. Then 
\[
\tilde f(x)=2x^{5^0+1}+x^{5^1+1}+x^{5^2+1}+2x^{5^3+1}+2x^{5^4+1}.
\]
Since
\[
S(\tilde f,1)=\sum_{x\in\Bbb F_5}e_1\bigl(\tilde f(x)\bigr)=\sum_{x\in\Bbb F_5}e_1(3x^2)=-g_5,
\]
we have $t_1(\tilde f)=-1$ and $l_1(\tilde f)=0$. By \eqref{5.5},
\[
t_{2^a}(f)=-1\qquad\text{for}\ a>0.
\]
Let $m=2^a13^bm^*$, where $s>0$ and $(m^*,2\cdot 5\cdot 13)=1$. Then by \eqref{4.7.0},
\[
t_m(f)=t_{2^a}(f)=-1.
\]
\end{exmp}

\begin{exmp}\label{E5.3}
\rm
Let
\[
f(x)=x^{3^0+1}+2x^{3^1+1}+2x^{3^2+1}+2x^{3^3+1}+x^{3^4+1}\in\Bbb F_3[x].
\]
Then
\[
f^*(x)=x+2x^3+2x^{3^2}+2x^{3^3}+2x^{3^4}+2x^{3^5}+2x^{3^6}+2x^{3^7}+x^{3^8}.
\]
The splitting field of $f^*$ over $\Bbb F_3$ is $\Bbb F_{3^{24}}$ and
\[
l_{2^a3^bm^*}(f)=
\begin{cases}
0&\text{if}\ a=0,\cr
1&\text{if}\ a=1,\ b=0,\cr
2&\text{if}\ a=1,\ b\ge 1,\cr
3&\text{if}\ a=2,\ b=0,\cr
4&\text{if}\ a=2,\ b\ge 1,\cr
7&\text{if}\ a\ge 3,\ b=0,\cr
8&\text{if}\ a\ge 3,\ b\ge 1,
\end{cases}
\]
where $(m^*,2\cdot 3)=1$. (See Table~\ref{Tb1} in the appendix.) Since
\[
S(f,1)=\sum_{x\in\Bbb F_3}e_1(2x^2)=-g_3,
\]
we have $t_1(f)=-1$. By Theorem~\ref{T4.3},
\[
t_{m^*}(f)=t_1(f)=-1\qquad\text{for all $m^*>0$ with $(m^*,\, 2\cdot 3)=1$.}
\]
Also,
\[
\tilde f(x)=2x^{3^0+1}+2x^{3^1+1}+x^{3^2+1}+2x^{3^3+1}+2x^{3^4+1}.
\]
Since $S(\tilde f,1)=\sum_{x\in\Bbb F_3}e_1(0)=3$, we have $t_1(\tilde f)=1$ and $l_1(\tilde f)=1$. By Theorem~\ref{T5.1},
\[
t_{2^a}(f)=-(-1)^a=(-1)^{a+1}\qquad\text{for all}\ a>0.
\]
By Theorem~\ref{T4.3} again,
\[
t_{2^am^*}(f)=(-1)^{a+1}\Bigl(\frac {m^*}3\Bigr),\qquad a>0,\ (m^*,2\cdot 3)=1.
\]
\end{exmp}

\begin{exmp}\label{E5.4}
\rm
Let
\[
f(x)=5x^{7^0+1}+6x^{7^1+1}+x^{7^2+1}\in\Bbb F_7[x].
\]
Using Mathematica, we find that the splitting field of $f^*$ over $\Bbb F_7$ is $\Bbb F_{7^{56}}$ and 
\[
l_{2^a7^bm^*}(f)=
\begin{cases}
0&\text{if}\ a=0,\cr
1&\text{if}\ 1\le a\le 2,\ b=0,\cr
2&\text{if}\ 1\le a\le 2,\ b\ge 1,\cr
3&\text{if}\ a\ge 3,\ b=0,\cr
4&\text{if}\ a\ge 3,\ b\ge 1,
\end{cases}
\]
where $(m^*,2\cdot 7)=1$. Since
\[
S(f,1)=\sum_{x\in\Bbb F_7}e_1(5x^2)=-g_7,
\]
we have $t_1(f)=-1$. By Theorem~\ref{T4.3},
\[
t_{m^*}(f)=t_1(f)=-1\qquad\text{for all $m^*>0$ with $(m^*,2\cdot 7)=1$.}
\]
Choose a nonsquare $\beta=3\in\Bbb F_7$. Then
\[
\tilde f(x)=x^{7^0+1}+3x^{7^1+1}+3x^{7^2+1}.
\]
Since $S(\tilde f,1)=\sum_{x\in\Bbb F_7}e_1(0)=7$, we have $t_1(\tilde f)=1$ and $l_1(\tilde f)=1$. By Theorem~\ref{T5.1},
\[
t_{2^a}(f)=-1\qquad\text{for all}\ a>0.
\]
By Theorem~\ref{T4.3} again,
\[
t_{2^am^*}(f)=-\Bigl(\frac{m^*}7\Bigr),\qquad a>0,\ (m^*,2\cdot 7)=1.
\]
\end{exmp}


\section{From $S(f,n)$ to $S(f,pn)$}

In this section, we consider the relative formula of $S(f,pn)$ in terms of $S(f,n)$. This seems to be a difficult situation; the congruence method in
the previous two sections does not work here. We are able to express $S(f,pn)$ in terms of $S(f,n)$ only under a very restrictive condition on $f$.

Let $b\in\Bbb F_{p^n}$ such that $\text{Tr}_n(b)\ne 0$. By the Artin-Schreier theorem (cf. \cite[Ch. VI, Theorem 6.4]{Lan93}), $\Bbb F_{p^{pn}}=\Bbb F_{p^n}(\epsilon)$,
where $\epsilon^p=\epsilon+b$, and the roots of $x^p-x-b$ are $\epsilon+j$, $j\in\Bbb F_p$. For every integer $t\ge 0$, we have
\begin{equation}\label{6.1}
\begin{split}
\text{Tr}_{pn/n}(\epsilon^t)\,&=\sum_{j\in\Bbb F_p}(\epsilon+j)^t=\sum_{j\in\Bbb F_p}\sum_{s=0}^t\binom ts j^s\epsilon^{t-s}\cr
&=\sum_{s=0}^t\binom ts\epsilon^{t-s}\sum_{j\in\Bbb F_p}j^s\cr
&=-\sum_{i>0}\binom t{i(p-1)}\epsilon^{t-i(p-1)}.
\end{split}
\end{equation}

\begin{lem}\label{L6.1}
Let $u,v$ be integers such that $0\le u,v\le p-1$. Then
\begin{equation}\label{6.2}
\text{\rm Tr}_{pn/n}(\epsilon^{u+v})=
\begin{cases}
0&\text{if}\ u+v\ne p-1,\ 2(p-1),\cr
-1&\text{if}\ u+v= p-1,\ 2(p-1).
\end{cases}
\end{equation}
If $\alpha$ is a positive integer, then
\begin{equation}\label{6.3}
\begin{split}
\text{\rm Tr}_{pn/n}(\epsilon^{u+vp^\alpha})=\,&-\binom v{u+v-(p-1)}(b^{p^0}+\cdots+b^{p^{\alpha-1}})^{u+v-(p-1)}\cr
&+\begin{cases}
0&\text{if}\ u+v\ne 2(p-1),\cr
-1&\text{if}\ u+v= 2(p-1).
\end{cases}
\end{split}
\end{equation}
\end{lem}

\begin{proof}
Write $u+v=u'+v'p$, where $0\le u',v'\le p-1$. By \eqref{6.1}, we have
\[
\begin{split}
\text{Tr}_{pn/n}(\epsilon^{u+v})\,&=-\binom{u+v}{p-1}\epsilon^{u+v-(p-1)}-\binom{u+v}{2(p-1)}
\epsilon^{u+v-2(p-1)}\cr
&=-\binom{u'+v'p}{p-1}\epsilon^{u+v-(p-1)}-\binom{u'+v'p}{(p-2)+1p}
\epsilon^{u+v-2(p-1)}\cr
&=-\binom{u'}{p-1}\epsilon^{u+v-(p-1)}-\binom{u'}{p-2}\binom{v'}1\epsilon^{u+v-2(p-1)}.
\end{split}
\]
Equation \eqref{6.2} follows immediately from the above.

Now we prove \eqref{6.3}. By \eqref{6.1},
\[
\begin{split}
&\text{Tr}_{pn/n}(\epsilon^{u+vp^\alpha})\cr
=\,&-\sum_{\substack{0\le s,t\le p-1\cr 0<s+tp^\alpha\equiv 0\kern-2mm \pmod{(p-1)}}}
\binom{u+vp^\alpha}{s+tp^\alpha}\epsilon^{u+vp^\alpha-(s+tp^\alpha)}\cr
=\,&-\sum_{\substack{0\le s\le u\cr 0\le t\le v\cr 0<s+t\equiv 0\kern-2mm \pmod{(p-1)}}}
\binom us \binom vt\epsilon^{u-s+(v-t)p^\alpha}\cr
=\,&-\sum_{\substack{0\le s\le u\cr 0\le t\le v\cr u+v>s+t\equiv u+v\kern-2mm \pmod{(p-1)}}}
\binom us \binom vt\epsilon^{s+tp^\alpha}\qquad (s\mapsto u-s,\ t\mapsto v-t).
\end{split}
\]
Since $\epsilon^p=\epsilon+b$, we have $\epsilon^{p^\alpha}=\epsilon+b^{p^0}+\cdots+b^{p^{\alpha-1}}$. Thus,
\[
\begin{split}
&\text{Tr}_{pn/n}(\epsilon^{u+vp^\alpha})\cr
=\,&-\sum_{\substack{0\le s\le u\cr 0\le t\le v\cr u+v>s+t\equiv u+v\kern-2mm \pmod{(p-1)}}}
\binom us \binom vt\epsilon^s(\epsilon+b^{p^0}+\cdots+b^{p^{\alpha-1}})^t\cr
=\,&-\sum_{\substack{0\le s\le u\cr 0\le t\le v\cr u+v>s+t\equiv u+v\kern-2mm \pmod{(p-1)}}}
\binom us \binom vt\sum_{\tau=0}^t\binom t\tau(b^{p^0}+\cdots+b^{p^{\alpha-1}})^{t-\tau}\epsilon^{s+\tau}.
\end{split}
\]
In the above sum, $s+\tau\le s+t\le u+v-(p-1)\le p-1$. Since $\text{Tr}_{pn/n}(\epsilon^{u+vp^\alpha})\in\Bbb F_{p^n}$, in the above, we only have to sum the terms with $s+\tau=0$. Therefore,
\[
\begin{split}
&\text{Tr}_{pn/n}(\epsilon^{u+vp^\alpha})\cr
=\,&-\sum_{\substack{0\le t\le v\cr u+v>t\equiv u+v\kern-2mm \pmod{(p-1)}}}
\binom vt(b^{p^0}+\cdots+b^{p^{\alpha-1}})^t\cr
=\,&-\binom v{u+v-(p-1)}(b^{p^0}+\cdots+b^{p^{\alpha-1}})^{u+v-(p-1)}+
\begin{cases}
0&\text{if}\ u+v\ne 2(p-1),\cr
-1&\text{if}\ u+v=2(p-1).
\end{cases}
\end{split}
\]
This completes the proof of the lemma.
\end{proof}

\begin{thm}\label{T6.2}
Let $f$ be as in \eqref{2.1.0}. Assume $\nu_p(n)<\min\{\nu_p(\alpha_i):1\le i\le k\}$.
Then,
\begin{equation}\label{6.3}
S(f,pn)=p^{\frac 12(p-3)(n+l_n(f))}|S(f,n)|^2\,\overline{S(f,n)}.
\end{equation}
Moreover, 
\begin{align}
\label{6.4}
l_{pn}(f)\,&=p\,l_n(f),\\
\label{6.5}
t_{pn}(f)\,&=t_n(f).
\end{align}
\end{thm}

\begin{proof}
Let $\nu_p(n)=\nu$ and write $n=p^\nu n'$, $p\nmid n'$. Choose $b\in\Bbb F_{p^{p^\nu}}$ such that $\text{Tr}_{p^\nu}(b)\ne 0$. Then $\text{Tr}_n(b)=n'\text{Tr}_{p^\nu}(b)\ne 0$.
Since $\nu_p(\alpha_i)>\nu$, we have
\[
b^{p^0}+\cdots+b^{p^{\alpha_i-1}}= \text{Tr}_{\alpha_i}(b)= 0.
\]
(If $\alpha_i=0$, $b^{p^0}+\cdots+b^{p^{\alpha_i-1}}$ is an empty sum.) By Lemma~\ref{L6.1},
for all $1\le i\le k$ and all $0\le u,v\le p-1$,
\[
\text{Tr}_{pn/n}(\epsilon^{u+vp^{\alpha_i}})=
\begin{cases}
0&\text{if}\ u+v\ne p-1,\ 2(p-1),\cr
-1&\text{if}\ u+v= p-1,\ 2(p-1).
\end{cases}
\]
Let $x=x_0\epsilon^0+\cdots+x_{p-1}\epsilon^{p-1}\in\Bbb F_{p^{pn}}$, where $x_u\in\Bbb F_{p^n}$,
$0\le u\le p-1$. Then 
\[
\begin{split}
&\text{Tr}_{pn}\bigl(f(x)\bigr)\cr
=\,&\text{Tr}_{pn}\Bigl(\sum_{i=1}^ka_i(x_0\epsilon^0+\cdots+x_{p-1}\epsilon^{p-1})^{1+p^{\alpha_i}}\Bigr)\cr
=\,&\text{Tr}_{pn}\Bigl(\sum_{i=1}^ka_i\sum_{0\le u,v\le p-1} x_ux_v^{p^{\alpha_i}}\epsilon^{u+vp^{\alpha_i}}\Bigr)\cr
=\,&\text{Tr}_n\Bigl[\sum_{i=1}^ka_i\sum_{0\le u,v\le p-1} x_ux_v^{p^{\alpha_i}}\text{Tr}_{pn/n}(\epsilon^{u+vp^{\alpha_i}})\Bigr]\cr
=\,&-\text{Tr}_n\Bigl[\sum_{i=1}^ka_i\sum_{\substack{0\le u,v\le p-1\cr
u+v=p-1,\,2(p-1)}}x_ux_v^{p^{\alpha_i}}\Bigr]\cr
=\,&-\text{Tr}_n\Bigl[\sum_{i=1}^ka_i\Bigl(x_{p-1}^{1+p^{\alpha_i}}+x_0x_{p-1}^{p^{\alpha_i}}
+x_{p-1}x_0^{p^{\alpha_i}}+x_{\frac 12(p-1)}^{1+p^{\alpha_i}}\cr
&+\sum_{u=1}^{\frac 12(p-3)}
(x_ux_{p-1-u}^{p^{\alpha_i}}+x_{p-1-u}x_u^{p^{\alpha_i}})\Bigr)\Bigr]\cr
=\,&-\text{Tr}_n\Bigl[\sum_{i=1}^ka_i\Bigl((x_0+x_{p-1})^{1+p^{\alpha_i}}-
x_0^{1+p^{\alpha_i}}+x_{\frac 12(p-1)}^{1+p^{\alpha_i}}\Bigr)\cr
&+\sum_{u=1}^{\frac 12(p-3)}\sum_{i=1}^ka_i(x_ux_{p-1-u}^{p^{\alpha_i}}+x_{p-1-u}x_u^{p^{\alpha_i}})\Bigr]\cr
=\,&-\text{Tr}_n\Bigl[f(x_0+x_{p-1})-f(x_0)+f(x_{\frac 12(p-1)})+\sum_{u=1}^{\frac 12(p-3)}
x_uf^*(x_{p-1-u})^{p^{-\alpha}}\Bigr].
\end{split}
\]
Therefore,
\[
\begin{split}
S(f,pn)
=\,&\sum_{(x_0,\dots,x_{p-1})\in\Bbb F_{p^n}^p}e_n
\Bigl[-f(x_0+x_{p-1})+f(x_0)-f(x_{\frac 12(p-1)})\cr
&-\sum_{u=1}^{\frac 12(p-3)}x_uf^*(x_{p-1-u})^{p^{-\alpha}}\Bigr]\cr
=\,&\overline{S(f,n)}^2S(f,n)\prod_{u=1}^{\frac 12(p-3)}\Bigl(\sum_{x_u,x_{p-1-u}\in\Bbb F_{p^n}}e_n\bigl(-x_u f^*(x_{p-1-u})^{p^{-\alpha}}\bigr)\Bigr)\cr
=\,&|S(f,n)|^2\,\overline{S(f,n)}\prod_{u=1}^{\frac 12(p-3)} p^{l_n(f)+n}\cr
=\,&p^{\frac 12(p-3)(n+l_n(f))}|S(f,n)|^2\,\overline{S(f,n)}.
\end{split}
\]
So, \eqref{6.3} is proved.

Equation~\eqref{6.3} gives
\[
t_{pn}(f)g_p^{pn-l_{pn}(f)}p^{l_{pn}(f)}=p^{\frac 12(p-3)(n+l_n(f))}(p^{\frac 12(n+l_n(f))})^2
t_n(f)\overline{g_p}^{n-l_n(f)}p^{l_n(f)},
\]
i.e.,
\[
t_{pn}(f) i^{\frac 14(p-1)^2(pn-l_{pn}(f))}p^{\frac 12(pn+l_{pn}(f))}
=t_n(f) i^{-\frac 14(p-1)^2(n-l_n(f))}p^{\frac 12p(n+l_n(f))}.
\]
Thus, $l_{pn}(f)=p\,l_n(f)$ and
\[
\begin{split}
t_{pn}(f)\,&=t_n(f) i^{-\frac 14(p-1)^2[n-l_n(f)+(pn-l_{pn}(f))]}\cr
&=t_n(f)\, i^{-\frac 14(p-1)^2(p+1)(n-l_n(f))}\cr
&=t_n(f).
\end{split}
\]
\end{proof}

\begin{cor}\label{C6.3}
In Theorem~\ref{T6.2}, assume $0\le s\le \min\{\nu_p(\alpha_i):1\le i\le k\}-\nu_p(n)$. Then
\begin{gather*}
l_{p^sn}(f)=p^sl_n(f),\\
t_{p^sn}(f)=t_n(f).
\end{gather*}
\end{cor}

\begin{proof}
Apply \eqref{6.4} and \eqref{6.5} repeatedly.
\end{proof}


\section{When $\nu_2(\alpha_1)=\cdots=\nu_2(\alpha_k)$}

\begin{lem}\label{L7.1}
Let $\alpha_1,\dots,\alpha_k\ge 0$ be integers. Then
\[
\text{\rm gcd}(p^{\alpha_1}+1,\dots,p^{\alpha_k}+1)>2\Leftrightarrow
\nu_2(\alpha_1)=\cdots=\nu_2(\alpha_k)<\infty.
\]
When $\nu_2(\alpha_1)=\cdots=\nu_2(\alpha_k)<\infty$,
\[
\text{\rm gcd}(p^{\alpha_1}+1,\dots,p^{\alpha_k}+1)=p^{\text{\rm gcd}(\alpha_1,\dots,\alpha_k)}+1.
\]
\end{lem}

\begin{proof}
It suffices to prove the lemma with $k=2$.

($\Leftarrow$) 
Since $\frac{\alpha_i}{(\alpha_1,\alpha_2)}$, $i=1,2$, are odd, $p^{(\alpha_1,\alpha_2)}+1\mid
p^{\alpha_i}+1$ for $i=1,2$. Thus, $p^{(\alpha_1,\alpha_2)}+1\mid (p^{\alpha_1}+1,
p^{\alpha_2}+1)$.

On the other hand,
\[
\begin{array}{rcl}
\displaystyle \frac 12(p^{\alpha_1}+1,p^{\alpha_2}+1) &\mid& \displaystyle \frac 12\Bigl(\frac 12(p^{2\alpha_1}-1),
\frac 12(p^{2\alpha_2}-1)\Bigr) \vrule height 0pt width0pt depth 12pt\cr
&=&\displaystyle \frac 14 (p^{2(\alpha_1,\alpha_2)}-1)\vrule height 0pt width0pt depth 12pt\cr
&=&\displaystyle \frac{p^{(\alpha_1,\alpha_2)}+1}2\cdot\frac{p^{(\alpha_1,\alpha_2)}-1}2.
\end{array}
\]
Since $(\frac{p^{\alpha_1}+1}2, \frac{p^{(\alpha_1,\alpha_2)}-1}2)=1$, we must have 
$\frac 12(p^{\alpha_1}+1,p^{\alpha_2}+1)\mid \frac 12(p^{(\alpha_1,\alpha_2)}+1)$,
i.e., $(p^{\alpha_1}+1,p^{\alpha_2}+1)\mid p^{(\alpha_1,\alpha_2)}+1$. Hence, 
$(p^{\alpha_1}+1,p^{\alpha_2}+1)= p^{(\alpha_1,\alpha_2)}+1$.

($\Rightarrow$) Clearly, $\alpha_i>0$ for all $1\le i\le k$. Assume to the contrary that $\nu_2(\alpha_1)>\nu_2(\alpha_2)$. Write $\alpha_1=2^i\alpha_1'$ and $\alpha_2=2^j\alpha_2'$, where $i>j$, $\alpha_1'$ and $\alpha_2'$ are odd. Then we have 
\[
\begin{array}{rcl}
(p^{\alpha_1}+1,p^{\alpha_2}+1) &\mid & (p^{2^i\alpha_1'\alpha_2'}+1, p^{2\alpha_2}-1)\cr
&\mid & (p^{2^i\alpha_1'\alpha_2'}+1, p^{2^i\alpha_1'\alpha_2'}-1)\cr
&=& 2,
\end{array}
\]
which is a contradiction.
\end{proof}

\begin{lem}\label{L7.2}
Let $\alpha,\beta\ge 0$ be integers. Then
\[
(p^\alpha+1,p^\beta-1)=
\begin{cases}
p^{(\alpha,\beta)}+1&\text{if}\ \nu_2(\beta)>\nu_2(\alpha),\cr
2&\text{if}\ \nu_2(\beta)\le \nu_2(\alpha).
\end{cases}
\]
\end{lem}

\begin{proof}
If any of $\alpha$ and $\beta$ is $0$, the conclusion is obvious. So, assume $\alpha,\beta>0$.

First assume $\nu_2(\beta)>\nu_2(\alpha)$. Since $\frac \alpha{(\alpha,\beta)}$ is odd,
$p^{(\alpha,\beta)}+1\mid p^\alpha+1$; since $\frac\beta{(\alpha,\beta)}$ is even,
$p^{(\alpha,\beta)}+1\mid p^\beta-1$. So, $p^{(\alpha,\beta)}+1\mid(p^\alpha+1,p^\beta-1)$.
Note that 
\[
\begin{array}{rcl}
\displaystyle \frac 12(p^\alpha+1,p^\beta-1) &\mid & \displaystyle \frac 12(p^{2\alpha}-1,p^\beta-1)\vrule height 0pt width0pt depth 12pt\cr
&=& \displaystyle \frac 12(p^{(2\alpha,\beta)}-1)\vrule height 0pt width0pt depth 12pt\cr
&=& \displaystyle \frac 12(p^{2(\alpha,\beta)}-1)\vrule height 0pt width0pt depth 12pt\cr
&=& \displaystyle \frac 12(p^{(\alpha,\beta)}-1)(p^{(\alpha,\beta)}+1).
\end{array}
\]
Since $(\frac 12(p^\alpha+1),\frac 12(p^{(\alpha,\beta)}-1))=1$, we must have
\begin{equation}\label{7.1}
\frac 12(p^\alpha+1,p^\beta-1)\mid p^{(\alpha,\beta)}+1.
\end{equation}
For each $x\in\Bbb Z$ and odd integer $k>0$, we have $\nu_2(1+x^k)=\nu_2(1+x)$. By this fact,
$\nu_2(p^\alpha+1)=\nu_2(p^{(\alpha,\beta)}+1)$ and $\nu_2(p^\beta-1)\ge\nu_2(p^\alpha+1)$.
So, $\nu_2(p^\alpha+1,p^\beta-1)=\nu_2(p^{(\alpha,\beta)}+1)$. Thus, \eqref{7.1} gives 
$(p^\alpha+1,p^\beta-1)\mid p^{(\alpha,\beta)}+1$. So, we have proved that $(p^\alpha+1,p^\beta-1)=p^{(\alpha,\beta)}+1$.

Now assume $\nu_2(\beta)\le\nu_2(\alpha)$. We have
\[
\begin{array}{rcl}
\displaystyle \Bigl(\frac{p^\alpha+1}2\cdot\frac{p^\alpha-1}2,\ \frac{p^\beta-1}2\Bigr) &\mid&
\displaystyle \frac 12 (p^{2\alpha}-1,p^\beta-1)\vrule height 0pt width0pt depth 12pt\cr
&=&\displaystyle \frac 12 (p^{(2\alpha,\beta)}-1)\vrule height 0pt width0pt depth 12pt\cr
&=&\displaystyle \frac 12 (p^{(\alpha,\beta)}-1)\vrule height 0pt width0pt depth 12pt\cr
&=&\displaystyle  \Bigl(\frac{p^\alpha-1}2,\ \frac{p^\beta-1}2\Bigr).
\end{array}
\]
Since $(\frac{p^\alpha+1}2,\, \frac{p^\alpha-1}2)=1$, we must have  
$(\frac{p^\alpha+1}2,\, \frac{p^\beta-1}2)=1$, i.e.,
$(p^\alpha+1,p^\beta-1)=2$.
\end{proof}

\begin{thm}\label{T7.3}
Let $f$ be as in \eqref{2.1.0}. Assume that $\nu_2(\alpha_1)=\cdots=\nu_2(\alpha_k)=\nu$ and
that $\nu_2(n)>\nu$. Then $2^{\nu+1}\mid l_n(f)$ and 
\[
t_n(f)=
(-1)^{(\frac 14(p-1)^22^\nu+1)\frac{n-l_n(f)}{2^{\nu+1}}}=
\begin{cases}
(-1)^{(\frac 14(p-1)^2+1)\frac{n-l_n(f)}{2^{\nu+1}}}&\text{if}\ \nu=0,\cr
(-1)^{\frac{n-l_n(f)}{2^{\nu+1}}}&\text{if}\ \nu>0.
\end{cases}
\]
\end{thm}

\begin{proof}
By Lemmas~\ref{L7.1} and \ref{7.2}, we have
\begin{equation}\label{7.2}
\begin{split}
\text{gcd}(p^{\alpha_1}+1,\dots,p^{\alpha_k}+1, p^n-1)\,&=\bigl(
p^{\text{gcd}(\alpha_1,\dots,\alpha_k)}+1,p^n-1\bigr)\cr
&= p^{\text{gcd}(\alpha_1,\dots,\alpha_k,n)}+1\cr
&\equiv 0\pmod{p^{2^\nu}+1}.
\end{split}
\end{equation}
Let $q=p^{2^\nu}+1$. Then $2^{\nu+1}$ is the multiplicative order
of $p$ modulo $q$, i.e., $o_q(p)=2^{\nu+1}$. Since $2^{\nu+1}\mid n$, we have $q\mid p^n-1$.

We first show that $2^{\nu+1}\mid l_n(f)$. In fact, we show that $\{x\in\Bbb F_{p^n}:f^*(x)=0\}$
is a vector space over $\Bbb F_{p^{2^{\nu+1}}}$. Let $x\in\Bbb F_{p^n}$ such that $f^*(x)=0$
and let $y\in\Bbb F_{p^{2^{\nu+1}}}$. We want to show $f^*(yx)=0$. By \eqref{2.3},
\begin{equation}\label{7.3}
f^*(yx)^{p^{-\alpha}}=\sum_{i=1}^k\bigl( a_iy^{p^{\alpha_i}}x^{p^{\alpha_i}}+
a_i^{p^{-\alpha_i}}y^{p^{-\alpha_i}}x^{p^{-\alpha_i}}\bigr).
\end{equation}
We claim that $y^{p^{\pm\alpha_i}}=y^{p^\alpha}$ for all $1\le i\le k$. By \eqref{7.2},
$p^{\alpha_i}\equiv -1\pmod q$; hence $p^{\pm\alpha_i}\equiv -1\pmod q$. Thus, $p^\alpha\equiv
-1\equiv p^{\pm\alpha_i}\pmod q$, i.e., $p^{\alpha\mp\alpha_i}\equiv 1\pmod q$ for all 
$1\le i\le k$. Since $o_q(p)=2^{\nu+1}$, $\alpha\mp\alpha_i\equiv 0\pmod{2^{\nu+1}}$. So,
$y^{p^{\alpha\mp\alpha_i}}=y$, i.e., $y^{p^\alpha}=y^{p^{\pm\alpha_i}}$.
Now by \eqref{7.3},
\[
f^*(yx)^{p^{-\alpha}}=y^{p^\alpha}\sum_{i=1}^k(a_ix^{p^{\alpha_i}}+a_i^{p^{-\alpha_i}}
x^{p^{-\alpha_i}})=y^{p^\alpha}f^*(x)^{p^{-\alpha}}=0.
\]

Choose $z\in \Bbb F_{p^n}^*$ such that $o(z)=q$. Since $p^{\alpha_i}+1\equiv 0\pmod q$, $1\le i\le k$, we have
\[
f(yx)=f(x)\qquad\text{for all}\ y\in\langle z\rangle.
\]
Thus,
\begin{equation}\label{7.4}
\begin{split}
t_n(f) g_p^{n-l_n(f)}p^{l_n(f)}\,&=S(f,n)\cr
&=1+\sum_{x\in\Bbb F_{p^n}^*}e_n(f(x))\cr
&=1+q\sum_{x\in\Bbb F_{p^n}^*/\langle z\rangle}e_n(f(x))\cr
&\equiv 1\pmod q.
\end{split}
\end{equation}
In the above, $p^{l_n(f)}\equiv 1\pmod q$ since $o_q(p)=2^{\nu+1}\mid l_n(f)$. Also,
\[
g_p^{2^{\nu+1}}=\Bigl[i^{\frac 14(p-1)^2}p^{\frac 12}\Bigr]^{2^{\nu+1}}=(-1)^{\frac 14(p-1)^22^\nu}p^{2^\nu}\equiv (-1)^{\frac 14(p-1)^2
2^\nu+1}\pmod q;
\]
hence,
\[
g_p^{n-l_n(f)}=g_p^{2^{\nu+1}\frac{n-l_n(f)}{2^{\nu+1}}}\equiv (-1)^{(\frac 14(p-1)^22^\nu+1)
\frac{n-l_n(f)}{2^{\nu+1}}}\pmod q.
\]
Now \eqref{7.4} gives
\[
t_n(f)=(-1)^{(\frac 14(p-1)^22^\nu+1)\frac{n-l_n(f)}{2^{\nu+1}}}.
\]
\end{proof}

\begin{cor}\label{C7.4}
Assume $p\equiv -1\pmod 4$, $\alpha_1,\dots,\alpha_k$ are all odd and $n$ is even. Then
$t_n(f)=1$.
\end{cor}

\begin{proof}
This immediate from Theorem~\ref{T7.3}.
\end{proof}

\begin{exmp}\label{E7.5}
\rm
Let $f(x)=3x^{5^1+1}+x^{5^3+1}\in\Bbb F_5[x]$. Then $f^*(x)=x+3x^{5^2}+3x^{5^4}+x^{5^6}$.
The splitting field of $f^*$ over $\Bbb F_5$ is $\Bbb F_{5^{20}}$ and
\[
l_{2^a5^bm^*}(f)=
\begin{cases}
0&\text{if}\ a\le 1,\cr
2&\text{if}\ a\ge 2,\ b=0,\cr
6&\text{if}\ a\ge 2,\ b=1,
\end{cases}
\]
where $(m^*,2\cdot 5)=1$. By Theorem~\ref{T7.3}, $t_n(f)=-1$ for all even $n$.
\end{exmp}

Theorem~\ref{T7.3} provides a quick proof for the explicit evaluation of the sum $S(ax^{p^\alpha+1},n)$.

\begin{cor}\label{C7.6}
Let $a\in \Bbb F_{p^n}^*$ and let $\alpha\ge 0$.
\begin{itemize}
\item[(i)]
If $\nu_2(n)\le\nu_2(\alpha)$,
\[
S(ax^{p^\alpha+1},n)=\eta(a)(-1)^{n-1}i^{\frac 14(p-1)^2n}p^{\frac 12 n}.
\]

\item[(ii)]
If
$\nu_2(n)=\nu_2(\alpha)+1$,
\[
S(ax^{p^\alpha+1},n)=
\begin{cases}
p^{\frac 12[n+(2\alpha,n)]}&\text{if}\ a^{\frac{(p^\alpha-1)(p^n-1)}{p^{(2\alpha,n)}-1}}=-1,\cr
-p^{\frac 12n}&\text{otherwise}.
\end{cases}
\]

\item[(iii)]
If $\nu_2(n)>\nu_2(\alpha)+1$,
\[
S(ax^{p^\alpha+1},n)=
\begin{cases}
-p^{\frac 12[n+(2\alpha,n)]}&\text{if}\ a^{\frac{(p^\alpha-1)(p^n-1)}{p^{(2\alpha,n)}-1}}=1,\cr
p^{\frac 12n}&\text{otherwise}.
\end{cases}
\]
\end{itemize}
\end{cor}

\begin{proof}
We first show that 
\begin{equation}\label{7.5}
l_n(f)=
\begin{cases}
(2\alpha,n)&\text{if}\ a^{\frac{(p^\alpha-1)(p^n-1)}{p^{(2\alpha,n)}-1}}=(-1)^{\frac{p^n-1}{p^{(2\alpha,n)}-1}},\cr
0&\text{otherwise}.
\end{cases}
\end{equation}
Note that $f^*(x)=a^{p^\alpha}x^{p^{2\alpha}}+ax=a^{p^\alpha}x(x^{p^{2\alpha}-1}+a^{1-p^\alpha})$. Thus, if $f^*(x)=0$ has a solution in $\Bbb F_{p^n}^*$,
the number of solutions is $(p^{2\alpha}-1,p^n-1)=p^{(2\alpha,n)}-1$. So, 
\[
l_n(f)=
\begin{cases}
(2\alpha,n)&\text{if $f^*(x)=0$ has a solution in $\Bbb F_{p^n}^*$},\cr
0&\text{otherwise}.
\end{cases}
\]
Observe that $f^*(x)=0$ has a solution in $\Bbb F_{p^n}^*$ if and only if $-a^{p^\alpha-1}=x^{p^{2\alpha}-1}$ for some $x\in\Bbb F_{p^n}^*$; 
the latter holds if and only if 
\[
\Bigl(-a^{p^\alpha-1}\Bigr)^{\frac{p^n-1}{(p^{2\alpha}-1,p^n-1)}}=1,
\]
i.e.,
\[
a^{\frac {(p^\alpha-1)(p^n-1)}{p^{(2\alpha,n)}-1}}=(-1)^{\frac{p^n-1}{p^{(2\alpha,n)}-1}}.
\]
So, \eqref{7.5} is proved.

(i) Since $\nu_2(n)\le \nu_2(\alpha)$, by Lemma~\ref{L7.2}, $(p^\alpha+1,p^n-1)=2$. Therefore, $x\mapsto x^{p^\alpha+1}$ is a 2-to-1 map from $\Bbb F_{p^n}^*$ to $(\Bbb F_{p^n}^*)^2$. Hence,
\[
\begin{split}
&S(ax^{p^\alpha+1},n)\cr
=\,& 1+\sum_{x\in\Bbb F_{p^n}^*}e_n(ax^{p^\alpha+1})=1+2\sum_{x\in(\Bbb F_{p^n}^*)^2} e_n(ax)\cr
=\,& 1+\sum_{x\in\Bbb F_{p^n}^*}e_n(ax^2)=\sum_{x\in\Bbb F_{p^n}}e_n(ax^2)=\eta(a)\sum_{x\in\Bbb F_{p^n}}e_n(x^2)\cr
=\,&\eta(a)(-1)^{n-1}g_p^n\qquad\text{(by the Davenport-Hasse theorem \cite{Dav35}, \cite[\S5.2]{Lid97})}\cr
=\,&\eta(a)(-1)^{n-1} i^{\frac 14(p-1)^2n}p^{\frac 12 n}.
\end{split}
\]

(ii) Since $\nu_2(n)=\nu_2((2\alpha,n))$, $\frac {p^n-1}{p^{(2\alpha,n)}-1}$ is odd. By \eqref{7.5},
\[
l_n(f)=
\begin{cases}
(2\alpha,n)&\text{if}\ a^{\frac {(p^\alpha-1)(p^n-1)}{p^{(2\alpha,n)}-1}}=-1,\cr
0&\text{otherwise}.
\end{cases}
\]
By Theorem~\ref{T7.3},
\[
t_n(f)=(-1)^{(\frac 14(p-1)^2\alpha+1)\frac{n-l_n(f)}{2^{\nu_2(\alpha)+1}}}.
\]
Thus,
\[
\begin{split}
S(ax^{p^\alpha+1},n)\,&=t_n(f) g_p^{n-l_n(f)}p^{l_n(f)}\cr
&=\begin{cases}
p^{\frac 12[n+(2\alpha,n)]}&\text{if}\ a^{\frac {(p^\alpha-1)(p^n-1)}{p^{(2\alpha,n)}-1}}=-1,\cr
-p^{\frac 12n}&\text{otherwise}.
\end{cases}
\end{split}
\]
(In the above, we used the fact that $\nu_2(n-(2\alpha,n))>\nu_2(n)$ and $\nu_2(2\alpha+n)>\nu_2(n)$.)

(iii) In this case, $\frac {p^n-1}{p^{(2\alpha,n)}-1}$ is even. By \eqref{7.5},
\[
l_n(f)=
\begin{cases}
(2\alpha,n)&\text{if}\ a^{\frac {(p^\alpha-1)(p^n-1)}{p^{(2\alpha,n)}-1}}=1,\cr
0&\text{otherwise}.
\end{cases}
\]
By Theorem~\ref{T7.3},
\[
t_n(f)=(-1)^{(\frac 14(p-1)^2\alpha+1)\frac{n-l_n(f)}{2^{\nu_2(\alpha)+1}}}.
\]
The conclusion follows the same way as in (ii). (Note that in this case, $\nu_2(n-(2\alpha,n))=\nu_2(\alpha)+1$.)
\end{proof}


\section*{Appendix. The Nullity $l_m(f)$ with $p^n=3$, $\alpha\le 4$ and $p^n=5$, $\alpha\le 3$}

Let $f(x)=\sum_{i=0}^ka_ix^{p^i+1}\in\Bbb F_{p^n}[x]$ with $a_k\in\Bbb F_p^*$ and let $0<m\equiv 0\pmod n$.
Since $\text{Tr}_m(a_k^{-1}f)=a_k^{-1}\text{Tr}_m(f)$, we have $l_m(f)=l_m(a_k^{-1}f)$, where $a_k^{-1}f$ is monic. So, when computing $l_m(f)$ with $a_k\in\Bbb F_p^*$, we may assume $a_k=1$.

This appendix contains two tables. Table~\ref{Tb1} gives the values of $l_m(f)$ with $p^n=3$ and $\alpha\le 4$;
table~\ref{Tb2} gives the values of $l_m(f)$ with $p^n=5$ and $\alpha\le 3$. In both tables, the left column contains
the coefficients $a_0,\dots,a_k$ of $f(x)=\sum_{i=0}^ka_ix^{p^i+1}$ with $a_k=1$. The middle column is $s$, where
$\Bbb F_{p^s}$ is the splitting field of $f^*(x)$. The right column lists all pairs $(m,l_m(f))$ such that
$m\mid s$. The values of $l_m(f)$ for arbitrary $m$ follows from \eqref{2.6}, cf. Example~\ref{E2.2}.

\begin{table}
\caption{Values of $l_m(f)$ with $p^n=3$, $\alpha\le 4$}\label{Tb1}
\vskip-5mm
\[
\begin{tabular}{l|c|l}
\hline
$a_0,\dots,a_k$ & $s$ & $(m,l_m(f)),\ m\mid s$ \\ \hline
1\ & 1&
(1,0)\ \\
0\ 1\ & 4&
(1,0)\ (2,0)\ (4,2)\ \\
1\ 1\ & 6&
(1,0)\ (2,1)\ (3,0)\ (6,2)\ \\
2\ 1\ & 3&
(1,1)\ (3,2)\ \\
0\ 0\ 1\ & 8&
(1,0)\ (2,0)\ (4,0)\ (8,4)\ \\
1\ 0\ 1\ & 12&
(1,0)\ (2,0)\ (3,0)\ (4,2)\ (6,0)\ (12,4)\ \\
2\ 0\ 1\ & 6&
(1,1)\ (2,2)\ (3,2)\ (6,4)\ \\
0\ 1\ 1\ & 18&
(1,0)\ (2,1)\ (3,0)\ (6,3)\ (9,0)\ (18,4)\ \\
1\ 1\ 1\ & 12&
(1,1)\ (2,1)\ (3,2)\ (4,3)\ (6,2)\ (12,4)\ \\
2\ 1\ 1\ & 5&
(1,0)\ (5,4)\ \\
0\ 2\ 1\ & 9&
(1,1)\ (3,3)\ (9,4)\ \\
1\ 2\ 1\ & 12&
(1,0)\ (2,1)\ (3,0)\ (4,3)\ (6,2)\ (12,4)\ \\
2\ 2\ 1\ & 10&
(1,0)\ (2,0)\ (5,0)\ (10,4)\ \\
0\ 0\ 0\ 1\ & 12&
(1,0)\ (2,0)\ (3,0)\ (4,2)\ (6,0)\ (12,6)\ \\
1\ 0\ 0\ 1\ & 18&
(1,0)\ (2,1)\ (3,0)\ (6,3)\ (9,0)\ (18,6)\ \\
2\ 0\ 0\ 1\ & 9&
(1,1)\ (3,3)\ (9,6)\ \\
0\ 1\ 0\ 1\ & 8&
(1,0)\ (2,0)\ (4,2)\ (8,6)\ \\
1\ 1\ 0\ 1\ & 30&
(1,1)\ (2,1)\ (3,2)\ (5,1)\ (6,2)\ (10,5)\ (15,2)\ (30,6)\ \\
2\ 1\ 0\ 1\ & 30&
(1,0)\ (2,1)\ (3,0)\ (5,4)\ (6,2)\ (10,5)\ (15,4)\ (30,6)\ \\
0\ 2\ 0\ 1\ & 12&
(1,1)\ (2,2)\ (3,2)\ (4,4)\ (6,4)\ (12,6)\ \\
1\ 2\ 0\ 1\ & 13&
(1,0)\ (13,6)\ \\
2\ 2\ 0\ 1\ & 26&
(1,0)\ (2,0)\ (13,0)\ (26,6)\ \\
0\ 0\ 1\ 1\ & 30&
(1,0)\ (2,1)\ (3,0)\ (5,0)\ (6,2)\ (10,5)\ (15,0)\ (30,6)\ \\
1\ 0\ 1\ 1\ & 12&
(1,1)\ (2,1)\ (3,2)\ (4,3)\ (6,2)\ (12,6)\ \\
2\ 0\ 1\ 1\ & 28&
(1,0)\ (2,0)\ (4,0)\ (7,0)\ (14,0)\ (28,6)\ \\
0\ 1\ 1\ 1\ & 24&
(1,1)\ (2,1)\ (3,2)\ (4,1)\ (6,2)\ (8,5)\ (12,2)\ (24,6)\ \\
1\ 1\ 1\ 1\ & 36&
(1,0)\ (2,1)\ (3,0)\ (4,3)\ (6,3)\ (9,0)\ (12,5)\ (18,4)\ (36,6)\ \\
2\ 1\ 1\ 1\ & 7&
(1,0)\ (7,6)\ \\
0\ 2\ 1\ 1\ & 13&
(1,0)\ (13,6)\ \\
1\ 2\ 1\ 1\ & 20&
(1,0)\ (2,0)\ (4,2)\ (5,4)\ (10,4)\ (20,6)\ \\
2\ 2\ 1\ 1\ & 18&
(1,1)\ (2,2)\ (3,3)\ (6,5)\ (9,4)\ (18,6)\ \\
0\ 0\ 2\ 1\ & 15&
(1,1)\ (3,2)\ (5,5)\ (15,6)\ \\
1\ 0\ 2\ 1\ & 28&
(1,0)\ (2,0)\ (4,0)\ (7,0)\ (14,0)\ (28,6)\ \\
2\ 0\ 2\ 1\ & 12&
(1,0)\ (2,1)\ (3,0)\ (4,3)\ (6,2)\ (12,6)\ \\
0\ 1\ 2\ 1\ & 24&
(1,0)\ (2,1)\ (3,0)\ (4,1)\ (6,2)\ (8,5)\ (12,2)\ (24,6)\ \\
1\ 1\ 2\ 1\ & 14&
(1,0)\ (2,0)\ (7,0)\ (14,6)\ \\
2\ 1\ 2\ 1\ & 36&
(1,1)\ (2,1)\ (3,3)\ (4,3)\ (6,3)\ (9,4)\ (12,5)\ (18,4)\ (36,6)\ \\
0\ 2\ 2\ 1\ & 26&
(1,0)\ (2,0)\ (13,0)\ (26,6)\ \\
1\ 2\ 2\ 1\ & 18&
(1,1)\ (2,2)\ (3,2)\ (6,5)\ (9,2)\ (18,6)\ \\
2\ 2\ 2\ 1\ & 20&
(1,0)\ (2,0)\ (4,2)\ (5,0)\ (10,4)\ (20,6)\ \\ 
\hline
\end{tabular}
\]
\end{table}

\addtocounter{table}{-1}

\begin{table}
\caption{Continued}
\vskip-5mm
\[
\begin{tabular}{l|c|l}
\hline
$a_0,\dots,a_k$ & $s$ & $(m,l_m(f)),\ m\mid s$ \\ \hline
0\ 0\ 0\ 0\ 1\ & 16&
(1,0)\ (2,0)\ (4,0)\ (8,0)\ (16,8)\ \\
1\ 0\ 0\ 0\ 1\ & 24&
(1,0)\ (2,0)\ (3,0)\ (4,0)\ (6,0)\ (8,4)\ (12,0)\ (24,8)\ \\
2\ 0\ 0\ 0\ 1\ & 12&
(1,1)\ (2,2)\ (3,2)\ (4,4)\ (6,4)\ (12,8)\ \\
0\ 1\ 0\ 0\ 1\ & 90&
(1,0)\ (2,1)\ (3,0)\ (5,0)\ (6,3)\ (9,0)\ (10,5)\ (15,0)\ (18,4)\ (30,7)\ (45,
    0)\ (90,8)\ \\
1\ 1\ 0\ 0\ 1\ & 84&
(1,1)\ (2,1)\ (3,2)\ (4,1)\ (6,2)\ (7,1)\ (12,2)\ (14,1)\ (21,2)\ (28,7)\ (42,
    2)\ (84,8)\ \\
2\ 1\ 0\ 0\ 1\ & 52&
(1,0)\ (2,0)\ (4,2)\ (13,0)\ (26,6)\ (52,8)\ \\
0\ 2\ 0\ 0\ 1\ & 45&
(1,1)\ (3,3)\ (5,5)\ (9,4)\ (15,7)\ (45,8)\ \\
1\ 2\ 0\ 0\ 1\ & 84&
(1,0)\ (2,1)\ (3,0)\ (4,1)\ (6,2)\ (7,0)\ (12,2)\ (14,1)\ (21,0)\ (28,7)\ (42,
    2)\ (84,8)\ \\
2\ 2\ 0\ 0\ 1\ & 52&
(1,0)\ (2,0)\ (4,2)\ (13,6)\ (26,6)\ (52,8)\ \\
0\ 0\ 1\ 0\ 1\ & 36&
(1,0)\ (2,0)\ (3,0)\ (4,2)\ (6,0)\ (9,0)\ (12,6)\ (18,0)\ (36,8)\ \\
1\ 0\ 1\ 0\ 1\ & 24&
(1,1)\ (2,2)\ (3,2)\ (4,2)\ (6,4)\ (8,6)\ (12,4)\ (24,8)\ \\
2\ 0\ 1\ 0\ 1\ & 10&
(1,0)\ (2,0)\ (5,4)\ (10,8)\ \\
0\ 1\ 1\ 0\ 1\ & 36&
(1,1)\ (2,1)\ (3,3)\ (4,3)\ (6,3)\ (9,6)\ (12,5)\ (18,6)\ (36,8)\ \\
1\ 1\ 1\ 0\ 1\ & 41&
(1,0)\ (41,8)\ \\
2\ 1\ 1\ 0\ 1\ & 42&
(1,0)\ (2,1)\ (3,0)\ (6,2)\ (7,6)\ (14,7)\ (21,6)\ (42,8)\ \\
0\ 2\ 1\ 0\ 1\ & 36&
(1,0)\ (2,1)\ (3,0)\ (4,3)\ (6,3)\ (9,0)\ (12,5)\ (18,6)\ (36,8)\ \\
1\ 2\ 1\ 0\ 1\ & 82&
(1,0)\ (2,0)\ (41,0)\ (82,8)\ \\
2\ 2\ 1\ 0\ 1\ & 42&
(1,1)\ (2,1)\ (3,2)\ (6,2)\ (7,1)\ (14,7)\ (21,2)\ (42,8)\ \\
0\ 0\ 2\ 0\ 1\ & 18&
(1,1)\ (2,2)\ (3,3)\ (6,6)\ (9,4)\ (18,8)\ \\
1\ 0\ 2\ 0\ 1\ & 24&
(1,0)\ (2,0)\ (3,0)\ (4,2)\ (6,0)\ (8,6)\ (12,4)\ (24,8)\ \\
2\ 0\ 2\ 0\ 1\ & 20&
(1,0)\ (2,0)\ (4,0)\ (5,0)\ (10,0)\ (20,8)\ \\
0\ 1\ 2\ 0\ 1\ & 41&
(1,0)\ (41,8)\ \\
1\ 1\ 2\ 0\ 1\ & 60&
(1,0)\ (2,1)\ (3,0)\ (4,3)\ (5,4)\ (6,2)\ (10,5)\ (12,4)\ (15,4)\ (20,7)\ (30,
    6)\ (60,8)\ \\
2\ 1\ 2\ 0\ 1\ & 78&
(1,1)\ (2,1)\ (3,2)\ (6,2)\ (13,1)\ (26,7)\ (39,2)\ (78,8)\ \\
0\ 2\ 2\ 0\ 1\ & 82&
(1,0)\ (2,0)\ (41,0)\ (82,8)\ \\
1\ 2\ 2\ 0\ 1\ & 60&
(1,1)\ (2,1)\ (3,2)\ (4,3)\ (5,1)\ (6,2)\ (10,5)\ (12,4)\ (15,2)\ (20,7)\ (30,
    6)\ (60,8)\ \\
2\ 2\ 2\ 0\ 1\ & 78&
(1,0)\ (2,1)\ (3,0)\ (6,2)\ (13,6)\ (26,7)\ (39,6)\ (78,8)\ \\
0\ 0\ 0\ 1\ 1\ & 42&
(1,0)\ (2,1)\ (3,0)\ (6,2)\ (7,0)\ (14,7)\ (21,0)\ (42,8)\ \\
1\ 0\ 0\ 1\ 1\ & 78&
(1,1)\ (2,1)\ (3,2)\ (6,2)\ (13,1)\ (26,7)\ (39,2)\ (78,8)\ \\
2\ 0\ 0\ 1\ 1\ & 60&
(1,0)\ (2,0)\ (3,0)\ (4,2)\ (5,4)\ (6,0)\ (10,4)\ (12,4)\ (15,4)\ (20,6)\ (30,
    4)\ (60,8)\ \\
0\ 1\ 0\ 1\ 1\ & 39&
(1,1)\ (3,2)\ (13,7)\ (39,8)\ \\
1\ 1\ 0\ 1\ 1\ & 72&
(1,0)\ (2,1)\ (3,0)\ (4,1)\ (6,3)\ (8,5)\ (9,0)\ (12,3)\ (18,4)\ (24,7)\ (36,
    4)\ (72,8)\ \\
2\ 1\ 0\ 1\ 1\ & 52&
(1,0)\ (2,0)\ (4,2)\ (13,6)\ (26,6)\ (52,8)\ \\
0\ 2\ 0\ 1\ 1\ & 80&
(1,0)\ (2,0)\ (4,0)\ (5,0)\ (8,0)\ (10,0)\ (16,0)\ (20,0)\ (40,0)\ (80,8)\ \\
1\ 2\ 0\ 1\ 1\ & 30&
(1,0)\ (2,0)\ (3,0)\ (5,0)\ (6,0)\ (10,4)\ (15,0)\ (30,8)\ \\
2\ 2\ 0\ 1\ 1\ & 36&
(1,1)\ (2,2)\ (3,3)\ (4,4)\ (6,5)\ (9,4)\ (12,7)\ (18,6)\ (36,8)\ \\
0\ 0\ 1\ 1\ 1\ & 12&
(1,1)\ (2,1)\ (3,2)\ (4,3)\ (6,2)\ (12,8)\ \\
1\ 0\ 1\ 1\ 1\ & 41&
(1,0)\ (41,8)\ \\
2\ 0\ 1\ 1\ 1\ & 78&
(1,0)\ (2,1)\ (3,0)\ (6,2)\ (13,0)\ (26,7)\ (39,0)\ (78,8)\ \\
0\ 1\ 1\ 1\ 1\ & 60&
(1,0)\ (2,1)\ (3,0)\ (4,3)\ (5,0)\ (6,2)\ (10,5)\ (12,4)\ (15,0)\ (20,7)\ (30,
    6)\ (60,8)\ \\
1\ 1\ 1\ 1\ 1\ & 40&
(1,0)\ (2,0)\ (4,0)\ (5,4)\ (8,4)\ (10,4)\ (20,4)\ (40,8)\ \\
\hline
\end{tabular}
\]
\end{table}

\addtocounter{table}{-1}

\begin{table}
\caption{Continued}
\vskip-5mm
\[
\begin{tabular}{l|c|l}
\hline
$a_0,\dots,a_k$ & $s$ & $(m,l_m(f)),\ m\mid s$ \\ \hline
2\ 1\ 1\ 1\ 1\ & 9&
(1,1)\ (3,3)\ (9,8)\ \\
0\ 2\ 1\ 1\ 1\ & 28&
(1,0)\ (2,0)\ (4,2)\ (7,0)\ (14,0)\ (28,8)\ \\
1\ 2\ 1\ 1\ 1\ & 18&
(1,1)\ (2,2)\ (3,2)\ (6,5)\ (9,2)\ (18,8)\ \\
2\ 2\ 1\ 1\ 1\ & 80&
(1,0)\ (2,0)\ (4,0)\ (5,0)\ (8,0)\ (10,0)\ (16,0)\ (20,0)\ (40,0)\ (80,8)\ \\
0\ 0\ 2\ 1\ 1\ & 41&
(1,0)\ (41,8)\ \\
1\ 0\ 2\ 1\ 1\ & 36&
(1,0)\ (2,1)\ (3,0)\ (4,3)\ (6,3)\ (9,0)\ (12,7)\ (18,4)\ (36,8)\ \\
2\ 0\ 2\ 1\ 1\ & 90&
(1,1)\ (2,1)\ (3,3)\ (5,1)\ (6,3)\ (9,4)\ (10,5)\ (15,3)\ (18,4)\ (30,7)\ (45,
    4)\ (90,8)\ \\
0\ 1\ 2\ 1\ 1\ & 82&
(1,0)\ (2,0)\ (41,0)\ (82,8)\ \\
1\ 1\ 2\ 1\ 1\ & 24&
(1,1)\ (2,1)\ (3,2)\ (4,3)\ (6,2)\ (8,7)\ (12,4)\ (24,8)\ \\
2\ 1\ 2\ 1\ 1\ & 84&
(1,0)\ (2,1)\ (3,0)\ (4,1)\ (6,2)\ (7,0)\ (12,2)\ (14,1)\ (21,0)\ (28,7)\ (42,
    2)\ (84,8)\ \\
0\ 2\ 2\ 1\ 1\ & 30&
(1,1)\ (2,2)\ (3,2)\ (5,5)\ (6,4)\ (10,6)\ (15,6)\ (30,8)\ \\
1\ 2\ 2\ 1\ 1\ & 28&
(1,0)\ (2,0)\ (4,2)\ (7,6)\ (14,6)\ (28,8)\ \\
2\ 2\ 2\ 1\ 1\ & 40&
(1,0)\ (2,0)\ (4,0)\ (5,0)\ (8,0)\ (10,0)\ (20,0)\ (40,8)\ \\
0\ 0\ 0\ 2\ 1\ & 21&
(1,1)\ (3,2)\ (7,7)\ (21,8)\ \\
1\ 0\ 0\ 2\ 1\ & 78&
(1,0)\ (2,1)\ (3,0)\ (6,2)\ (13,6)\ (26,7)\ (39,6)\ (78,8)\ \\
2\ 0\ 0\ 2\ 1\ & 60&
(1,0)\ (2,0)\ (3,0)\ (4,2)\ (5,0)\ (6,0)\ (10,4)\ (12,4)\ (15,0)\ (20,6)\ (30,
    4)\ (60,8)\ \\
0\ 1\ 0\ 2\ 1\ & 80&
(1,0)\ (2,0)\ (4,0)\ (5,0)\ (8,0)\ (10,0)\ (16,0)\ (20,0)\ (40,0)\ (80,8)\ \\
1\ 1\ 0\ 2\ 1\ & 15&
(1,0)\ (3,0)\ (5,4)\ (15,8)\ \\
2\ 1\ 0\ 2\ 1\ & 36&
(1,1)\ (2,2)\ (3,2)\ (4,4)\ (6,5)\ (9,2)\ (12,7)\ (18,6)\ (36,8)\ \\
0\ 2\ 0\ 2\ 1\ & 78&
(1,0)\ (2,1)\ (3,0)\ (6,2)\ (13,0)\ (26,7)\ (39,0)\ (78,8)\ \\
1\ 2\ 0\ 2\ 1\ & 72&
(1,1)\ (2,1)\ (3,3)\ (4,1)\ (6,3)\ (8,5)\ (9,4)\ (12,3)\ (18,4)\ (24,7)\ (36,
    4)\ (72,8)\ \\
2\ 2\ 0\ 2\ 1\ & 52&
(1,0)\ (2,0)\ (4,2)\ (13,0)\ (26,6)\ (52,8)\ \\
0\ 0\ 1\ 2\ 1\ & 12&
(1,0)\ (2,1)\ (3,0)\ (4,3)\ (6,2)\ (12,8)\ \\
1\ 0\ 1\ 2\ 1\ & 82&
(1,0)\ (2,0)\ (41,0)\ (82,8)\ \\
2\ 0\ 1\ 2\ 1\ & 39&
(1,1)\ (3,2)\ (13,7)\ (39,8)\ \\
0\ 1\ 1\ 2\ 1\ & 28&
(1,0)\ (2,0)\ (4,2)\ (7,0)\ (14,0)\ (28,8)\ \\
1\ 1\ 1\ 2\ 1\ & 18&
(1,1)\ (2,2)\ (3,3)\ (6,5)\ (9,6)\ (18,8)\ \\
2\ 1\ 1\ 2\ 1\ & 80&
(1,0)\ (2,0)\ (4,0)\ (5,0)\ (8,0)\ (10,0)\ (16,0)\ (20,0)\ (40,0)\ (80,8)\ \\
0\ 2\ 1\ 2\ 1\ & 60&
(1,1)\ (2,1)\ (3,2)\ (4,3)\ (5,5)\ (6,2)\ (10,5)\ (12,4)\ (15,6)\ (20,7)\ (30,
    6)\ (60,8)\ \\
1\ 2\ 1\ 2\ 1\ & 40&
(1,0)\ (2,0)\ (4,0)\ (5,0)\ (8,4)\ (10,4)\ (20,4)\ (40,8)\ \\
2\ 2\ 1\ 2\ 1\ & 18&
(1,0)\ (2,1)\ (3,0)\ (6,3)\ (9,0)\ (18,8)\ \\
0\ 0\ 2\ 2\ 1\ & 82&
(1,0)\ (2,0)\ (41,0)\ (82,8)\ \\
1\ 0\ 2\ 2\ 1\ & 36&
(1,1)\ (2,1)\ (3,3)\ (4,3)\ (6,3)\ (9,4)\ (12,7)\ (18,4)\ (36,8)\ \\
2\ 0\ 2\ 2\ 1\ & 90&
(1,0)\ (2,1)\ (3,0)\ (5,4)\ (6,3)\ (9,0)\ (10,5)\ (15,4)\ (18,4)\ (30,7)\ (45,
    4)\ (90,8)\ \\
0\ 1\ 2\ 2\ 1\ & 30&
(1,1)\ (2,2)\ (3,2)\ (5,1)\ (6,4)\ (10,6)\ (15,2)\ (30,8)\ \\
1\ 1\ 2\ 2\ 1\ & 28&
(1,0)\ (2,0)\ (4,2)\ (7,0)\ (14,6)\ (28,8)\ \\
2\ 1\ 2\ 2\ 1\ & 40&
(1,0)\ (2,0)\ (4,0)\ (5,0)\ (8,0)\ (10,0)\ (20,0)\ (40,8)\ \\
0\ 2\ 2\ 2\ 1\ & 41&
(1,0)\ (41,8)\ \\
1\ 2\ 2\ 2\ 1\ & 24&
(1,0)\ (2,1)\ (3,0)\ (4,3)\ (6,2)\ (8,7)\ (12,4)\ (24,8)\ \\
2\ 2\ 2\ 2\ 1\ & 84&
(1,1)\ (2,1)\ (3,2)\ (4,1)\ (6,2)\ (7,1)\ (12,2)\ (14,1)\ (21,2)\ (28,7)\ (42,
    2)\ (84,8)\ \\
\hline
\end{tabular}
\]
\end{table}


\begin{table}
\caption{Values of $l_m(f)$ with $p^n=5$, $\alpha\le 3$}\label{Tb2}
\vskip-5mm
\[
\begin{tabular}{l|c|l}
\hline
$a_0,\dots,a_k$ & $s$ & $(m,l_m(f)),\ m\mid s$ \\ \hline
1\ & 1&
(1,0)\ \\
0\ 1\ & 4&
(1,0)\ (2,0)\ (4,2)\ \\
1\ 1\ & 10&
(1,0)\ (2,1)\ (5,0)\ (10,2)\ \\
2\ 1\ & 6&
(1,0)\ (2,0)\ (3,0)\ (6,2)\ \\
3\ 1\ & 3&
(1,0)\ (3,2)\ \\
4\ 1\ & 5&
(1,1)\ (5,2)\ \\
0\ 0\ 1\ & 8&
(1,0)\ (2,0)\ (4,0)\ (8,4)\ \\
1\ 0\ 1\ & 20&
(1,0)\ (2,0)\ (4,2)\ (5,0)\ (10,0)\ (20,4)\ \\
2\ 0\ 1\ & 12&
(1,0)\ (2,0)\ (3,0)\ (4,0)\ (6,0)\ (12,4)\ \\
3\ 0\ 1\ & 6&
(1,0)\ (2,0)\ (3,2)\ (6,4)\ \\
4\ 0\ 1\ & 10&
(1,1)\ (2,2)\ (5,2)\ (10,4)\ \\
0\ 1\ 1\ & 30&
(1,0)\ (2,1)\ (3,0)\ (5,0)\ (6,3)\ (10,2)\ (15,0)\ (30,4)\ \\
1\ 1\ 1\ & 12&
(1,0)\ (2,0)\ (3,2)\ (4,2)\ (6,2)\ (12,4)\ \\
2\ 1\ 1\ & 13&
(1,0)\ (13,4)\ \\
3\ 1\ 1\ & 5&
(1,1)\ (5,4)\ \\
4\ 1\ 1\ & 24&
(1,0)\ (2,0)\ (3,0)\ (4,0)\ (6,0)\ (8,0)\ (12,0)\ (24,4)\ \\
0\ 2\ 1\ & 26&
(1,0)\ (2,0)\ (13,0)\ (26,4)\ \\
1\ 2\ 1\ & 20&
(1,0)\ (2,1)\ (4,3)\ (5,0)\ (10,2)\ (20,4)\ \\
2\ 2\ 1\ & 30&
(1,1)\ (2,1)\ (3,1)\ (5,2)\ (6,3)\ (10,2)\ (15,2)\ (30,4)\ \\
3\ 2\ 1\ & 13&
(1,0)\ (13,4)\ \\
4\ 2\ 1\ & 15&
(1,0)\ (3,2)\ (5,0)\ (15,4)\ \\
0\ 3\ 1\ & 13&
(1,0)\ (13,4)\ \\
1\ 3\ 1\ & 20&
(1,1)\ (2,1)\ (4,3)\ (5,2)\ (10,2)\ (20,4)\ \\
2\ 3\ 1\ & 30&
(1,0)\ (2,1)\ (3,2)\ (5,0)\ (6,3)\ (10,2)\ (15,2)\ (30,4)\ \\
3\ 3\ 1\ & 26&
(1,0)\ (2,0)\ (13,0)\ (26,4)\ \\
4\ 3\ 1\ & 30&
(1,0)\ (2,0)\ (3,0)\ (5,0)\ (6,2)\ (10,0)\ (15,0)\ (30,4)\ \\
0\ 4\ 1\ & 15&
(1,1)\ (3,3)\ (5,2)\ (15,4)\ \\
1\ 4\ 1\ & 12&
(1,0)\ (2,0)\ (3,0)\ (4,2)\ (6,2)\ (12,4)\ \\
2\ 4\ 1\ & 26&
(1,0)\ (2,0)\ (13,0)\ (26,4)\ \\
3\ 4\ 1\ & 10&
(1,0)\ (2,1)\ (5,0)\ (10,4)\ \\
4\ 4\ 1\ & 24&
(1,0)\ (2,0)\ (3,0)\ (4,0)\ (6,0)\ (8,0)\ (12,0)\ (24,4)\ \\
0\ 0\ 0\ 1\ & 12&
(1,0)\ (2,0)\ (3,0)\ (4,2)\ (6,0)\ (12,6)\ \\
1\ 0\ 0\ 1\ & 30&
(1,0)\ (2,1)\ (3,0)\ (5,0)\ (6,3)\ (10,2)\ (15,0)\ (30,6)\ \\
2\ 0\ 0\ 1\ & 18&
(1,0)\ (2,0)\ (3,0)\ (6,0)\ (9,0)\ (18,6)\ \\
3\ 0\ 0\ 1\ & 9&
(1,0)\ (3,0)\ (9,6)\ \\
4\ 0\ 0\ 1\ & 15&
(1,1)\ (3,3)\ (5,2)\ (15,6)\ \\
0\ 1\ 0\ 1\ & 8&
(1,0)\ (2,0)\ (4,2)\ (8,6)\ \\
1\ 1\ 0\ 1\ & 126&
(1,0)\ (2,0)\ (3,0)\ (6,0)\ (7,0)\ (9,0)\ (14,0)\ (18,0)\ (21,0)\ (42,0)\ (63,
    0)\ (126,6)\ \\
2\ 1\ 0\ 1\ & 30&
(1,0)\ (2,1)\ (3,2)\ (5,0)\ (6,3)\ (10,4)\ (15,2)\ (30,6)\ \\
3\ 1\ 0\ 1\ & 30&
(1,1)\ (2,1)\ (3,1)\ (5,4)\ (6,3)\ (10,4)\ (15,4)\ (30,6)\ \\
4\ 1\ 0\ 1\ & 63&
(1,0)\ (3,0)\ (7,0)\ (9,0)\ (21,0)\ (63,6)\ \\
0\ 2\ 0\ 1\ & 12&
(1,0)\ (2,0)\ (3,2)\ (4,2)\ (6,4)\ (12,6)\ \\
\hline
\end{tabular}
\]
\end{table}

\addtocounter{table}{-1}

\begin{table}
\caption{Continued}
\vskip-5mm
\[
\begin{tabular}{l|c|l}
\hline
$a_0,\dots,a_k$ & $s$ & $(m,l_m(f)),\ m\mid s$ \\ \hline
1\ 2\ 0\ 1\ & 21&
(1,0)\ (3,0)\ (7,0)\ (21,6)\ \\
2\ 2\ 0\ 1\ & 130&
(1,1)\ (2,1)\ (5,2)\ (10,2)\ (13,1)\ (26,5)\ (65,2)\ (130,6)\ \\
3\ 2\ 0\ 1\ & 130&
(1,0)\ (2,1)\ (5,0)\ (10,2)\ (13,4)\ (26,5)\ (65,4)\ (130,6)\ \\
4\ 2\ 0\ 1\ & 42&
(1,0)\ (2,0)\ (3,0)\ (6,0)\ (7,0)\ (14,0)\ (21,0)\ (42,6)\ \\
0\ 3\ 0\ 1\ & 20&
(1,0)\ (2,0)\ (4,2)\ (5,0)\ (10,0)\ (20,6)\ \\
1\ 3\ 0\ 1\ & 65&
(1,1)\ (5,2)\ (13,5)\ (65,6)\ \\
2\ 3\ 0\ 1\ & 24&
(1,0)\ (2,0)\ (3,0)\ (4,0)\ (6,2)\ (8,0)\ (12,2)\ (24,6)\ \\
3\ 3\ 0\ 1\ & 24&
(1,0)\ (2,0)\ (3,2)\ (4,0)\ (6,2)\ (8,0)\ (12,2)\ (24,6)\ \\
4\ 3\ 0\ 1\ & 130&
(1,0)\ (2,1)\ (5,0)\ (10,2)\ (13,0)\ (26,5)\ (65,0)\ (130,6)\ \\
0\ 4\ 0\ 1\ & 20&
(1,1)\ (2,2)\ (4,4)\ (5,2)\ (10,4)\ (20,6)\ \\
1\ 4\ 0\ 1\ & 78&
(1,0)\ (2,0)\ (3,2)\ (6,2)\ (13,0)\ (26,4)\ (39,2)\ (78,6)\ \\
2\ 4\ 0\ 1\ & 62&
(1,0)\ (2,0)\ (31,0)\ (62,6)\ \\
3\ 4\ 0\ 1\ & 31&
(1,0)\ (31,6)\ \\
4\ 4\ 0\ 1\ & 78&
(1,0)\ (2,0)\ (3,0)\ (6,2)\ (13,4)\ (26,4)\ (39,4)\ (78,6)\ \\
0\ 0\ 1\ 1\ & 50&
(1,0)\ (2,1)\ (5,0)\ (10,5)\ (25,0)\ (50,6)\ \\
1\ 0\ 1\ 1\ & 52&
(1,0)\ (2,0)\ (4,2)\ (13,4)\ (26,4)\ (52,6)\ \\
2\ 0\ 1\ 1\ & 12&
(1,0)\ (2,0)\ (3,2)\ (4,0)\ (6,2)\ (12,6)\ \\
3\ 0\ 1\ 1\ & 65&
(1,1)\ (5,2)\ (13,5)\ (65,6)\ \\
4\ 0\ 1\ 1\ & 78&
(1,0)\ (2,0)\ (3,0)\ (6,2)\ (13,4)\ (26,4)\ (39,4)\ (78,6)\ \\
0\ 1\ 1\ 1\ & 24&
(1,0)\ (2,0)\ (3,2)\ (4,0)\ (6,2)\ (8,4)\ (12,2)\ (24,6)\ \\
1\ 1\ 1\ 1\ & 60&
(1,0)\ (2,1)\ (3,0)\ (4,3)\ (5,0)\ (6,3)\ (10,2)\ (12,5)\ (15,0)\ (20,4)\ (30,
    4)\ (60,6)\ \\
2\ 1\ 1\ 1\ & 130&
(1,1)\ (2,1)\ (5,2)\ (10,2)\ (13,1)\ (26,5)\ (65,2)\ (130,6)\ \\
3\ 1\ 1\ 1\ & 7&
(1,0)\ (7,6)\ \\
4\ 1\ 1\ 1\ & 62&
(1,0)\ (2,0)\ (31,0)\ (62,6)\ \\
0\ 2\ 1\ 1\ & 63&
(1,0)\ (3,0)\ (7,0)\ (9,0)\ (21,0)\ (63,6)\ \\
1\ 2\ 1\ 1\ & 20&
(1,1)\ (2,1)\ (4,3)\ (5,4)\ (10,4)\ (20,6)\ \\
2\ 2\ 1\ 1\ & 120&
(1,0)\ (2,1)\ (3,0)\ (4,1)\ (5,0)\ (6,1)\ (8,1)\ (10,2)\ (12,1)\ (15,0)\ (20,
    2)\ (24,5)\\
&& (30,2)\ (40,2)\ (60,2)\ (120,6)\ \\
3\ 2\ 1\ 1\ & 30&
(1,0)\ (2,0)\ (3,2)\ (5,0)\ (6,4)\ (10,0)\ (15,4)\ (30,6)\ \\
4\ 2\ 1\ 1\ & 31&
(1,0)\ (31,6)\ \\
0\ 3\ 1\ 1\ & 30&
(1,1)\ (2,1)\ (3,1)\ (5,2)\ (6,3)\ (10,2)\ (15,2)\ (30,6)\ \\
1\ 3\ 1\ 1\ & 60&
(1,0)\ (2,0)\ (3,2)\ (4,2)\ (5,0)\ (6,2)\ (10,0)\ (12,4)\ (15,2)\ (20,4)\ (30,
    2)\ (60,6)\ \\
2\ 3\ 1\ 1\ & 42&
(1,0)\ (2,0)\ (3,0)\ (6,0)\ (7,0)\ (14,0)\ (21,0)\ (42,6)\ \\
3\ 3\ 1\ 1\ & 130&
(1,0)\ (2,1)\ (5,0)\ (10,2)\ (13,0)\ (26,5)\ (65,0)\ (130,6)\ \\
4\ 3\ 1\ 1\ & 124&
(1,0)\ (2,0)\ (4,0)\ (31,0)\ (62,0)\ (124,6)\ \\
0\ 4\ 1\ 1\ & 124&
(1,0)\ (2,0)\ (4,0)\ (31,0)\ (62,0)\ (124,6)\ \\
1\ 4\ 1\ 1\ & 24&
(1,0)\ (2,0)\ (3,0)\ (4,2)\ (6,0)\ (8,2)\ (12,2)\ (24,6)\ \\
2\ 4\ 1\ 1\ & 78&
(1,0)\ (2,0)\ (3,0)\ (6,2)\ (13,0)\ (26,4)\ (39,0)\ (78,6)\ \\
3\ 4\ 1\ 1\ & 124&
(1,0)\ (2,0)\ (4,0)\ (31,0)\ (62,0)\ (124,6)\ \\
4\ 4\ 1\ 1\ & 30&
(1,1)\ (2,2)\ (3,3)\ (5,2)\ (6,4)\ (10,4)\ (15,4)\ (30,6)\ \\
0\ 0\ 2\ 1\ & 24&
(1,0)\ (2,0)\ (3,2)\ (4,0)\ (6,2)\ (8,0)\ (12,2)\ (24,6)\ \\
\hline
\end{tabular}
\]
\end{table}

\addtocounter{table}{-1}

\begin{table}
\caption{Continued}
\vskip-5mm
\[
\begin{tabular}{l|c|l}
\hline
$a_0,\dots,a_k$ & $s$ & $(m,l_m(f)),\ m\mid s$ \\ \hline
1\ 0\ 2\ 1\ & 63&
(1,0)\ (3,0)\ (7,0)\ (9,0)\ (21,0)\ (63,6)\ \\
2\ 0\ 2\ 1\ & 60&
(1,1)\ (2,1)\ (3,1)\ (4,3)\ (5,2)\ (6,3)\ (10,2)\ (12,5)\ (15,2)\ (20,4)\ (30,
    4)\ (60,6)\ \\
3\ 0\ 2\ 1\ & 124&
(1,0)\ (2,0)\ (4,0)\ (31,0)\ (62,0)\ (124,6)\ \\
4\ 0\ 2\ 1\ & 60&
(1,0)\ (2,1)\ (3,0)\ (4,1)\ (5,0)\ (6,1)\ (10,2)\ (12,5)\ (15,0)\ (20,2)\ (30,
    2)\ (60,6)\ \\
0\ 1\ 2\ 1\ & 40&
(1,0)\ (2,1)\ (4,1)\ (5,0)\ (8,5)\ (10,2)\ (20,2)\ (40,6)\ \\
1\ 1\ 2\ 1\ & 120&
(1,1)\ (2,1)\ (3,1)\ (4,1)\ (5,2)\ (6,1)\ (8,1)\ (10,2)\ (12,1)\ (15,2)\ (20,
    2)\ (24,5)\\
&& (30,2)\ (40,2)\ (60,2)\ (120,6)\ \\
2\ 1\ 2\ 1\ & 52&
(1,0)\ (2,0)\ (4,2)\ (13,0)\ (26,4)\ (52,6)\ \\
3\ 1\ 2\ 1\ & 39&
(1,0)\ (3,2)\ (13,4)\ (39,6)\ \\
4\ 1\ 2\ 1\ & 30&
(1,0)\ (2,0)\ (3,0)\ (5,0)\ (6,2)\ (10,0)\ (15,0)\ (30,6)\ \\
0\ 2\ 2\ 1\ & 130&
(1,1)\ (2,1)\ (5,2)\ (10,2)\ (13,1)\ (26,5)\ (65,2)\ (130,6)\ \\
1\ 2\ 2\ 1\ & 30&
(1,0)\ (2,1)\ (3,2)\ (5,0)\ (6,5)\ (10,2)\ (15,2)\ (30,6)\ \\
2\ 2\ 2\ 1\ & 52&
(1,0)\ (2,0)\ (4,2)\ (13,4)\ (26,4)\ (52,6)\ \\
3\ 2\ 2\ 1\ & 62&
(1,0)\ (2,0)\ (31,0)\ (62,6)\ \\
4\ 2\ 2\ 1\ & 63&
(1,0)\ (3,0)\ (7,0)\ (9,0)\ (21,0)\ (63,6)\ \\
0\ 3\ 2\ 1\ & 63&
(1,0)\ (3,0)\ (7,0)\ (9,0)\ (21,0)\ (63,6)\ \\
1\ 3\ 2\ 1\ & 126&
(1,0)\ (2,0)\ (3,0)\ (6,0)\ (7,0)\ (9,0)\ (14,0)\ (18,0)\ (21,0)\ (42,0)\ (63,
    0)\ (126,6)\ \\
2\ 3\ 2\ 1\ & 20&
(1,0)\ (2,1)\ (4,3)\ (5,0)\ (10,2)\ (20,6)\ \\
3\ 3\ 2\ 1\ & 78&
(1,0)\ (2,0)\ (3,0)\ (6,2)\ (13,4)\ (26,4)\ (39,4)\ (78,6)\ \\
4\ 3\ 2\ 1\ & 15&
(1,1)\ (3,3)\ (5,4)\ (15,6)\ \\
0\ 4\ 2\ 1\ & 78&
(1,0)\ (2,0)\ (3,0)\ (6,2)\ (13,0)\ (26,4)\ (39,0)\ (78,6)\ \\
1\ 4\ 2\ 1\ & 62&
(1,0)\ (2,0)\ (31,0)\ (62,6)\ \\
2\ 4\ 2\ 1\ & 60&
(1,0)\ (2,0)\ (3,2)\ (4,2)\ (5,0)\ (6,2)\ (10,0)\ (12,4)\ (15,4)\ (20,2)\ (30,
    4)\ (60,6)\ \\
3\ 4\ 2\ 1\ & 10&
(1,1)\ (2,2)\ (5,2)\ (10,6)\ \\
4\ 4\ 2\ 1\ & 124&
(1,0)\ (2,0)\ (4,0)\ (31,0)\ (62,0)\ (124,6)\ \\
0\ 0\ 3\ 1\ & 24&
(1,0)\ (2,0)\ (3,0)\ (4,0)\ (6,2)\ (8,0)\ (12,2)\ (24,6)\ \\
1\ 0\ 3\ 1\ & 60&
(1,1)\ (2,1)\ (3,1)\ (4,1)\ (5,2)\ (6,1)\ (10,2)\ (12,5)\ (15,2)\ (20,2)\ (30,
    2)\ (60,6)\ \\
2\ 0\ 3\ 1\ & 124&
(1,0)\ (2,0)\ (4,0)\ (31,0)\ (62,0)\ (124,6)\ \\
3\ 0\ 3\ 1\ & 60&
(1,0)\ (2,1)\ (3,2)\ (4,3)\ (5,0)\ (6,3)\ (10,2)\ (12,5)\ (15,2)\ (20,4)\ (30,
    4)\ (60,6)\ \\
4\ 0\ 3\ 1\ & 126&
(1,0)\ (2,0)\ (3,0)\ (6,0)\ (7,0)\ (9,0)\ (14,0)\ (18,0)\ (21,0)\ (42,0)\ (63,
    0)\ (126,6)\ \\
0\ 1\ 3\ 1\ & 40&
(1,1)\ (2,1)\ (4,1)\ (5,2)\ (8,5)\ (10,2)\ (20,2)\ (40,6)\ \\
1\ 1\ 3\ 1\ & 15&
(1,0)\ (3,2)\ (5,0)\ (15,6)\ \\
2\ 1\ 3\ 1\ & 78&
(1,0)\ (2,0)\ (3,0)\ (6,2)\ (13,0)\ (26,4)\ (39,0)\ (78,6)\ \\
3\ 1\ 3\ 1\ & 52&
(1,0)\ (2,0)\ (4,2)\ (13,4)\ (26,4)\ (52,6)\ \\
4\ 1\ 3\ 1\ & 120&
(1,0)\ (2,1)\ (3,0)\ (4,1)\ (5,0)\ (6,1)\ (8,1)\ (10,2)\ (12,1)\ (15,0)\ (20,
    2)\ (24,5)\\
&& (30,2)\ (40,2)\ (60,2)\ (120,6)\ \\
0\ 2\ 3\ 1\ & 130&
(1,0)\ (2,1)\ (5,0)\ (10,2)\ (13,4)\ (26,5)\ (65,4)\ (130,6)\ \\
1\ 2\ 3\ 1\ & 126&
(1,0)\ (2,0)\ (3,0)\ (6,0)\ (7,0)\ (9,0)\ (14,0)\ (18,0)\ (21,0)\ (42,0)\ (63,
    0)\ (126,6)\ \\
2\ 2\ 3\ 1\ & 31&
(1,0)\ (31,6)\ \\
3\ 2\ 3\ 1\ & 52&
(1,0)\ (2,0)\ (4,2)\ (13,0)\ (26,4)\ (52,6)\ \\
4\ 2\ 3\ 1\ & 30&
(1,1)\ (2,1)\ (3,3)\ (5,2)\ (6,5)\ (10,2)\ (15,4)\ (30,6)\ \\
0\ 3\ 3\ 1\ & 126&
(1,0)\ (2,0)\ (3,0)\ (6,0)\ (7,0)\ (9,0)\ (14,0)\ (18,0)\ (21,0)\ (42,0)\ (63,
    0)\ (126,6)\ \\
\hline
\end{tabular}
\]
\end{table}

\addtocounter{table}{-1}

\begin{table}
\caption{Continued}
\vskip-5mm
\[
\begin{tabular}{l|c|l}
\hline
$a_0,\dots,a_k$ & $s$ & $(m,l_m(f)),\ m\mid s$ \\ \hline
1\ 3\ 3\ 1\ & 30&
(1,0)\ (2,1)\ (3,0)\ (5,0)\ (6,3)\ (10,4)\ (15,0)\ (30,6)\ \\
2\ 3\ 3\ 1\ & 78&
(1,0)\ (2,0)\ (3,2)\ (6,2)\ (13,0)\ (26,4)\ (39,2)\ (78,6)\ \\
3\ 3\ 3\ 1\ & 20&
(1,1)\ (2,1)\ (4,3)\ (5,2)\ (10,2)\ (20,6)\ \\
4\ 3\ 3\ 1\ & 63&
(1,0)\ (3,0)\ (7,0)\ (9,0)\ (21,0)\ (63,6)\ \\
0\ 4\ 3\ 1\ & 39&
(1,0)\ (3,2)\ (13,4)\ (39,6)\ \\
1\ 4\ 3\ 1\ & 124&
(1,0)\ (2,0)\ (4,0)\ (31,0)\ (62,0)\ (124,6)\ \\
2\ 4\ 3\ 1\ & 10&
(1,1)\ (2,2)\ (5,4)\ (10,6)\ \\
3\ 4\ 3\ 1\ & 60&
(1,0)\ (2,0)\ (3,0)\ (4,2)\ (5,0)\ (6,2)\ (10,0)\ (12,4)\ (15,0)\ (20,2)\ (30,
    4)\ (60,6)\ \\
4\ 4\ 3\ 1\ & 31&
(1,0)\ (31,6)\ \\
0\ 0\ 4\ 1\ & 25&
(1,1)\ (5,5)\ (25,6)\ \\
1\ 0\ 4\ 1\ & 78&
(1,0)\ (2,0)\ (3,2)\ (6,2)\ (13,0)\ (26,4)\ (39,2)\ (78,6)\ \\
2\ 0\ 4\ 1\ & 130&
(1,0)\ (2,1)\ (5,0)\ (10,2)\ (13,0)\ (26,5)\ (65,0)\ (130,6)\ \\
3\ 0\ 4\ 1\ & 12&
(1,0)\ (2,0)\ (3,0)\ (4,0)\ (6,2)\ (12,6)\ \\
4\ 0\ 4\ 1\ & 52&
(1,0)\ (2,0)\ (4,2)\ (13,0)\ (26,4)\ (52,6)\ \\
0\ 1\ 4\ 1\ & 24&
(1,0)\ (2,0)\ (3,0)\ (4,0)\ (6,2)\ (8,4)\ (12,2)\ (24,6)\ \\
1\ 1\ 4\ 1\ & 31&
(1,0)\ (31,6)\ \\
2\ 1\ 4\ 1\ & 14&
(1,0)\ (2,0)\ (7,0)\ (14,6)\ \\
3\ 1\ 4\ 1\ & 130&
(1,0)\ (2,1)\ (5,0)\ (10,2)\ (13,4)\ (26,5)\ (65,4)\ (130,6)\ \\
4\ 1\ 4\ 1\ & 60&
(1,1)\ (2,1)\ (3,3)\ (4,3)\ (5,2)\ (6,3)\ (10,2)\ (12,5)\ (15,4)\ (20,4)\ (30,
    4)\ (60,6)\ \\
0\ 2\ 4\ 1\ & 126&
(1,0)\ (2,0)\ (3,0)\ (6,0)\ (7,0)\ (9,0)\ (14,0)\ (18,0)\ (21,0)\ (42,0)\ (63,
    0)\ (126,6)\ \\
1\ 2\ 4\ 1\ & 62&
(1,0)\ (2,0)\ (31,0)\ (62,6)\ \\
2\ 2\ 4\ 1\ & 30&
(1,0)\ (2,0)\ (3,2)\ (5,0)\ (6,4)\ (10,0)\ (15,2)\ (30,6)\ \\
3\ 2\ 4\ 1\ & 120&
(1,1)\ (2,1)\ (3,1)\ (4,1)\ (5,2)\ (6,1)\ (8,1)\ (10,2)\ (12,1)\ (15,2)\ (20,
    2)\ (24,5)\\
&& (30,2)\ (40,2)\ (60,2)\ (120,6)\ \\
4\ 2\ 4\ 1\ & 20&
(1,0)\ (2,1)\ (4,3)\ (5,0)\ (10,4)\ (20,6)\ \\
0\ 3\ 4\ 1\ & 30&
(1,0)\ (2,1)\ (3,2)\ (5,0)\ (6,3)\ (10,2)\ (15,4)\ (30,6)\ \\
1\ 3\ 4\ 1\ & 124&
(1,0)\ (2,0)\ (4,0)\ (31,0)\ (62,0)\ (124,6)\ \\
2\ 3\ 4\ 1\ & 65&
(1,1)\ (5,2)\ (13,5)\ (65,6)\ \\
3\ 3\ 4\ 1\ & 21&
(1,0)\ (3,0)\ (7,0)\ (21,6)\ \\
4\ 3\ 4\ 1\ & 60&
(1,0)\ (2,0)\ (3,0)\ (4,2)\ (5,0)\ (6,2)\ (10,0)\ (12,4)\ (15,0)\ (20,4)\ (30,
    2)\ (60,6)\ \\
0\ 4\ 4\ 1\ & 124&
(1,0)\ (2,0)\ (4,0)\ (31,0)\ (62,0)\ (124,6)\ \\
1\ 4\ 4\ 1\ & 30&
(1,1)\ (2,2)\ (3,1)\ (5,2)\ (6,4)\ (10,4)\ (15,2)\ (30,6)\ \\
2\ 4\ 4\ 1\ & 124&
(1,0)\ (2,0)\ (4,0)\ (31,0)\ (62,0)\ (124,6)\ \\
3\ 4\ 4\ 1\ & 39&
(1,0)\ (3,2)\ (13,4)\ (39,6)\ \\
4\ 4\ 4\ 1\ & 24&
(1,0)\ (2,0)\ (3,0)\ (4,2)\ (6,0)\ (8,2)\ (12,2)\ (24,6)\ \\
\hline
\end{tabular}
\]
\end{table}


\end{document}